\documentclass[smallextended,referee,envcountsect,]{svjour3}
\smartqed 
\usepackage{amsmath}
\usepackage{amssymb}
\usepackage{enumitem}
\usepackage{tikz}
\usetikzlibrary{arrows,hobby}
\usetikzlibrary{patterns, patterns.meta}

\journalname{Preprint}

\begin{document}

\title{Existence of Solutions in Bi-level Stochastic Linear Programming with Integer Variables}

\titlerunning{Existence of Solutions for Bi-level Stochastic MILPs}

\author{Johanna Burtscheidt \and Matthias Claus}

\institute{Johanna Burtscheidt, Corresponding author \at
            University of Duisburg-Essen\\
            Essen, Germany\\
            johanna.burtscheidt@uni-due.de
        \and
            Matthias Claus \at
            University of Duisburg-Essen\\
            Essen, Germany\\
            matthias.claus@uni-due.de
}


\maketitle

\begin{abstract}
The addition of lower level integrality constraints to a bi-level linear program is known to result in significantly weaker analytical properties. Most notably, the upper level goal function in the optimistic setting lacks lower semicontinuity and the existence of an optimal solution cannot be guaranteed under standard assumptions.
In this paper, we study a setting where the right-hand side of the lower level constraint system is affected by the leader's choice as well as the realization of some random vector. Assuming that only the follower decides under complete information, we employ a convex risk measure to assess the upper level outcome. Confining the analysis to the cases where the lower level feasible set is finite, we provide sufficient conditions for Hölder continuity of the leader's risk functional and draw conclusions about the existence of optimal solutions. Finally, we examine qualitative stability with respect to perturbations of the underlying probability measure. Considering the topology of weak convergence, we prove joint continuity of the objective function with respect to both the leader's decision and the underlying probability measure.
\end{abstract}

\section{Introduction}
The present work deals with bi-level stochastic linear problems in which the realization of a random vector whose distribution does not depend on the upper level decision enters the lower level as an additional parameter. It is assumed that the leader must make his decision without knowing the realization of the random vector, while the follower decides under complete information. However, the variables controlled by the follower are integral, which entails that techniques from the non-integer case are not or only partially transferable to the present model. In addition, the feasible set of the follower is finite.
An application of this problem type is conceivable, for example, in the transport sector: The leader decides on the expansion of a road or rail network, i.e., how procured certain road sections should be built in order to protect them from flooding or other weather effects on the one hand and to minimize costs on the other hand? After the realization of the randomness (i.e. the water level or weather), the follower solves a shortest path problem: Are the extended route sections usable, or does he have to fall back on existing cab connections or bus shuttles, the use of which causes costs for the leader in the course of customer satisfaction?
Various types of price setting problems are listed in \cite{LaVi16}, where two types of network price problems use continuous decisions at the upper level and binary decisions at the lower level.
Applications that make use of a bi-level knapsack problem can be found in \cite[\S 5]{Fa06}, \cite{OePrSc10} and references there. In the case of \cite{OePrSc10} even bi-level knapsack problems with stochastic right hand side are considered: They are a stochastic extension of the bi-level knapsack problem, where the upper level decision has an uncertain effect on the knapsack capacity in the lower level. 

\smallskip
In the literature, bi-level problems with integer variables in the lower level are often studied in combination with a suitable solution approach. We first review some literature on bi-level (mixed-)integer linear programming with (mixed-)integer lower level:
In \cite{GuFl05}, Gümü\c{s} and Floudas give a classification of mixed-integer nonlinear bi-level programming problems and consider global solution strategies for all identified classes.
Dempe and Mefo Kue, Mefo Kue as well as Vicente, Savard and Judice, respectively, also give the above classification in \cite{DeMe17}, in \cite{Me17}, or in \cite{ViSaJu96}.
The latter consider equivalences between different classes of discrete linear bi-level programs and particular linear multi-level problems, and analyze some properties of the discrete linear bi-level program for different discretizations of the variable set.
In \cite{YuGaZeYo19}, Yue, Gao, Zeng and You propose a decomposition algorithm through column-and-constraint generation for mixed-integer linear bi-level problems based on a projection-based single-level formulation.
A mixed-integer linear bi-level problem solution algorithm was introduced in \cite{MoBa90} by Moore and Bard. 
Caramia and Mari propose two exact algorithms for integer linear bi-level problems based on a cutting plane method and a branch-and-cut algorithm in \cite{CaMa15}.
A branch-and-cut framework and an accompanying open source solver, MibS, for the pure integer bi-level linear program is provided by DeNegre in \cite{De11}.
In \cite{BrHaMa13}, Brotcorne, Hanafi, and Mansi describe a solution approach to bi-level knapsack problem in terms of an integer linear bi-level problem using dynamic programming and a reformulation into a linear model with integer variables.
Kara and Verter consider the problem of designing a road network for hazardous materials transportation and present a solution methodology for its underlying binary linear bi-level problem in \cite{KaVe04}.

We now summarize some existing work on problems where the leader's variables are continuous and follower's are integers. This part includes bi-level stochastic problems:
Bi-level problems with a continuous upper level and with a discrete or integer lower level are also considered in the aforementioned paper \cite{GuFl05} by Gümü\c{s} and Floudas.
Köppe, Queyranne and Ryan show in \cite{KoQuRy10} the existence of an algorithm for mixed-integer linear bi-level problems, where the leader’s variables are continuous and the follower solves an integer linear program, and one for integer linear bi-level problems.
For the same types of problems, Dempe and Mefo Kue establish algorithms in \cite{DeMe17} based on the optimal value reformulations of the bi-level problems and the use of branch-and-cut.
In her dissertation \cite{Me17}, Mefo Kue analyses bi-level problems with integer lower level and bi-level problems with a semidefinite programming problem in the lower level predominantly with regard to the existence of solutions and optimality conditions. She establishes two different solution approaches for bi-level problems with continuous upper level and integer lower level as well as an algorithm for solving the discrete linear bi-level problem where both levels are binary or bounded integer based on the optimal value reformulation.
Fanghänel's dissertation \cite{Fa06} deals, among other things, with (extended) solution sets, (weak) solution functions, optimality conditions, and a solution approach for weak local solutions of bi-level problems with a continuous upper level and a discrete lower level. She applies this knowledge to bi-level problems with discrete convex lower level problems, to the problem class of linear bi-level problems with a 0-1 knapsack problem at the lower level and to the assortment pricing problem.
Özalt\i n, Prokopyev, and Schaefer also consider a bi-level knapsack problem in \cite{OePrSc10}, but with stochastic right-hand side in the follower constraints. For a model with continuous leader variables and binary follower variables, they give a solvability result and formulate a solvability strategy for the case of integer leader variables based on a branch-and-backtrack algorithm and a branch-and-cut algorithm.
In \cite{ZhOe21}, Zhang and Özalt\i n study structural properties of a model with continuous upper level variables and integer lower level variables, where the right-hand side in the follower constraints is stochastic, too, e.g., they show that there exists a characterization of the bi-level integer problem value function only by bi-level minimal right-hand side vectors. They design a dynamic programming algorithm as well as a two-step approach to solve the value function reformulation of the original bi-level integer problem and finally apply it to a bi-level facility interdiction problem with stochastic resource constraints.

\smallskip
Under suitable assumptions, the optimal value function of an integer linear problem is known to be lower semicontinuous in the right-hand side of its constraint system. In an optimistic bi-level framework, a linear functional is optimized over the set of minimizers of said integer problem. Unfortunately, lower semicontinuity does not carry over to the resulting upper level goal function, which complicates the analysis. Assuming that only the leader decides nonanticipatorily, the presence of stochastic right-hand side uncertainty in the lower level gives rise to a family of random variables indexed by the leaders's decision. We study these random variables based on geometric insights from the underlying parametric problem.

Assessing the random upper level outcome by convex risk measueres, we establish sufficient conditions for Hölder continuity of the resulting objective function under weak assumptions. This allows us to formulate sufficient conditions for the existence of global minimizers despite the difficulties described above. Incomplete information or the need to compute efficiently can lead to optimization models that use an approximation to the true underlying distribution. This motivates the analysis of the behavior of optimal values and local optimal solution sets under perturbations of the underlying distribution. Finally, we establish a qualitative stability result for bi-level stochastic linear problems with integer variables that holds for all law-invariant convex risk measures.

\smallskip
The remainder of the paper is organized as follows. In Section~\ref{Sec_2}, we describe a parametric bi-level integer linear problem and characterize the domains of the feasible set mapping, as well as the goal functions of follower and leader. In the next section, we extend our model by introducing stochastic uncertainty and formulate a bi-level stochastic integer linear problem. From this, we consider the feasibility set induced by implicit constraints and determine the properties of the leader's goal function under various assumptions. Properties of the final risk-neutral or risk-averse problem, including an existence result, follow in Section~\ref{Sec_4}. Section~\ref{Sec_stability} contains the analysis of the bi-level stochastic integer linear bi-level problem in terms of joint continuity of the leader variable and the probability measure. The last section of this paper presents our conclusions.

\section{Structural properties of bi-level integer programs} \label{Sec_2}
Consider the optimistic bi-level linear program
\begin{equation}\label{BSILP}
    \min_x \left\{ c^\top x + \inf_y \left\{ q^\top y \; | \; y \in \Psi(Tx+z) \right\} \; | \; x \in X \right\},
\end{equation}
where $X \subseteq \mathbb{R}^n$ is a nonempty polyhedron and $z \in \mathbb{R}^s$ a parameter that enters the lower level optimal solution set mapping $\Psi: \mathbb{R}^s \rightrightarrows \mathbb{Z}^m$ defined by
$$
    \Psi(t) 
    := \underset{y}{\mathrm{Argmin}} \left\{ d^\top y \; | \; Wy \leq t, \, y \in Y \cap \mathbb{Z}^m \right\}.
$$
Throughout the paper we assume $Y \cap \mathbb{Z}^m \neq \emptyset$. In the set $Y \subseteq \mathbb{R}^m$ we summarize all constraints of the follower that do not depend on the parameter $t \in \mathbb{R}^s$, i.e. are independent of the decision of the leader $x$ and the parameter $z$. Typically, these are sign constraints or box constraints on the follower's variables $y$. For $W = 0$ or $q = 0$ or $d = 0$ or $q = d$ problem \eqref{BSILP} collapses to different types of single-level problems.
The following example shows that the reduced upper level goal function $\varphi: \mathbb{R}^s \to \overline{\mathbb{R}} := \mathbb{R} \cup \lbrace \pm \infty \rbrace$ with $\varphi(t)
    := \inf_y \left\{q^\top y \; | \; y \in \Psi(t) \right\}$ is neither lower nor upper semicontinuous nor convex in general:
\begin{example} \label{Ex_1}
Let $\varphi: \mathbb{R} \times \mathbb{R} \to \mathbb{R}$ be given by
\begin{align*}
    &\varphi(t_1, t_2)\\
    &= \inf_y \left\{y_1 + y_2 \; | \; (y_1, y_2) \in \underset{y'}{\mathrm{Argmin}} \left\{y_1' - y_2' \; | \; -y_1' \leq t_1, \, y_2' \leq t_2, \, y' \in \mathbb{Z}^2\right\}\right\}\\
    &= \lceil -t_1\rceil + \lfloor t_2\rfloor,
\end{align*}
then the restriction of $\varphi$ to the linear subspace $\lbrace (-t, t) \; | \; t \in \mathbb{R} \rbrace$ is neither upper nor lower semicontinuous nor convex at any integral point, cf. Figure~\ref{Fig_varphi(t,t)}.
\end{example}

We shall first establish some basic properties of $\varphi$:
\begin{lemma} \label{Lemma_phimeasurable2}
The extended real-valued function $\varphi$ is measurable.
\end{lemma}

{\it Proof} In view of the representation $\varphi(t) 
    = \inf_{y \in Y \cap \mathbb{Z}^m} q^\top y \cdot \chi_{\Psi(t)}(y)$ with
$$
    \chi_{\Psi(t)}(y) 
    = \begin{cases} 1, &\text{if} \; y \in \Psi(t) \\ \infty, &\text{else} \end{cases}
$$
and \cite[Th. 1.14]{Ru86}, it is sufficient to show that $\chi_{\Psi(t)}(y)$ is measurable for any fixed $y \in Y \cap \mathbb{Z}^m$. This is equivalent to measurability of
\begin{align*}
    &\lbrace t \; | \; y \in \Psi(t) \rbrace\\
    &= \lbrace t \; | \; Wy \leq t \rbrace \cap \bigcup_{y' \in Y \cap \mathbb{Z}^m} \left[ \lbrace t \; | \; d^\top y \leq d^\top y' \rbrace \; \cup \; \bigcup_{i=1, \ldots, s} \lbrace t \; | \; e_i^\top Wy' > e_i^\top t \rbrace \right],
\end{align*}
which follows from countability of $Y \cap \mathbb{Z}^m$.\qed

\smallskip
Let $\Phi: \mathbb{R}^s \rightrightarrows \mathbb{R}^m$ with $\Phi(t) 
    := \lbrace y \in Y \cap \mathbb{Z}^m \; | \; Wy \leq t \rbrace$ denote the lower level feasible set mapping and let $\psi: \mathbb{R}^s \to \overline{\mathbb{R}}$ with $\psi(t)
    := \inf \{d^\top y \; | \; y \in \Phi(t)\}$ be the lower level optimal value function. This function, unlike the non-integer case, cf. \cite[Th. 3.15]{De02}, is generally non-convex:

\medskip
\noindent \textit{Example~\ref{Ex_1} (continued) }
We get $\psi: \mathbb{R} \times \mathbb{R} \to \mathbb{R}$ with
\begin{align*}
    \psi(t_1, t_2) 
    = \inf_y \left\{y_1' - y_2' \; | \; -y_1' \leq t_1, \, y_2' \leq t_2, \, y' \in \mathbb{R}^2 \cap \mathbb{Z}^2\right\}
    = -\lfloor t_1\rfloor - \lfloor t_2\rfloor,
\end{align*}
then the restriction of $\psi$ to the linear subspace $\{(t, t) \; | \; t \in \mathbb{R}\}$ is non-convex at any integral point, cf. Figure~\ref{Fig_varphi(t,t)2}.
\begin{figure}[ht]
\centering
\begin{minipage}[t]{0.49\textwidth}
\centering
\vspace{-0.15 in}
\begin{tikzpicture}[scale=0.7]
\filldraw[black] (-1, -2)circle[radius=2pt]
(0, 0)circle[radius=2pt]
(1, 2)circle[radius=2pt]
(2, 4)circle[radius=2pt];
\foreach \x in {-1,1,2}
  \draw (\x,2pt) -- (\x,-2pt) node[below] {$\x$};
\foreach \y in {1,...,4}
  \draw (2pt,\y) -- (-2pt,\y) node[left] {$\y$};
\foreach \y in {-2,-1}
  \draw (2pt,\y) -- (-2pt,\y) node[right] {$\y$};
\draw[thick, ->] (-1.5, 0) -- (2.5, 0) node[below] {$t$};
\draw[thick, ->] (0, -2.5) -- (0, 4.5) node[left] {$\varphi(-t, t)$};
\draw[{(-)}] (-1, -1) -- (0, -1);
\draw[{(-)}] (0, 1) -- (1, 1);
\draw[{(-)}] (1, 3) -- (2, 3);
\end{tikzpicture}
\caption{Graph of the mapping \mbox{$t \mapsto \varphi(-t, t)$} of Example~\ref{Ex_1}}
\label{Fig_varphi(t,t)}
\vspace{-0.15 in}
\end{minipage}
\begin{minipage}[t]{0.49\textwidth}
\centering
\vspace{-0.15 in}
\begin{tikzpicture}[scale=0.7]
\filldraw[black] (-1, 2)circle[radius=2pt]
(0, 0)circle[radius=2pt]
(1, -2)circle[radius=2pt]
(2, -4)circle[radius=2pt];
\foreach \x in {-1,1,2}
  \draw (\x,2pt) -- (\x,-2pt) node[below] {$\x$};
\foreach \y in {1,2}
  \draw (2pt,\y) -- (-2pt,\y) node[right] {$\y$};
\foreach \y in {-4,...,-1}
  \draw (2pt,\y) -- (-2pt,\y) node[right] {$\y$};
\draw[thick, ->] (-1.5, 0) -- (2.5, 0) node[below] {$t$};
\draw[thick, ->] (0, -4.5) -- (0, 2.5) node[left] {$\psi(t, t)$};
\draw[{-)}] (-1, 2) -- (0, 2);
\draw[{-)}] (0, 0) -- (1, 0);
\draw[{-)}] (1, -2) -- (2, -2);
\draw[-] (2, -4) -- (2.5, -4);
\end{tikzpicture}
\caption{Graph of the mapping \mbox{$t \mapsto \psi(t, t)$} of Example~\ref{Ex_1}~(cont.)}
\label{Fig_varphi(t,t)2}
\vspace{-0.15 in}
\end{minipage}
\end{figure}

\medskip
We collect some properties of $\psi$ that can be formulated based on \cite[Prop. 1]{ZhOe21} and \cite[Th. 4.5.2.]{5M} including their proofs:
\begin{proposition}
Let $Y = \mathbb{R}^m$.
    \begin{enumerate}
        \item The function $\psi$ is nonincreasing on $\mathbb{R}^s$.
        \item The function $\psi$ is subadditive on $\mathrm{dom} \; \Phi$.
        \item If the matrix $W$ has only rational elements, then the function $\psi$ is lower semicontinuous on $\mathrm{dom} \; \Phi$.
    \end{enumerate}
\end{proposition}

Similar results are found in \cite[Prop. 3.1] {DeMe17} and \cite[Prop. 5.2.3.]{Me17} followed by the proofs.
Due to \cite[Rem. 1]{ZhOe21}, the following example shows that the reduced upper level goal function fails to be monotone and subadditive in general:
\begin{example} \label{Ex_subadd}
Let $\varphi(t_1, t_2) = \inf_y \left\{-3y_1 - 2y_2 - y_3 \; | \; y \in \Psi(t_1, t_2)\right\}$ with
\begin{align*}
    \left.\begin{aligned}
        \Psi(t_1, t_2)
        = \underset{y}{\mathrm{Argmin}} \big\{-y_1 - 2y_2 - 2y_3 \; | \; &y_1 + 2y_2 + 2y_3 \leq t_1,\\
        &2y_1 + 2y_2 + y_3 \leq t_2, \, y \in \mathbb{N}_0^3
    \end{aligned}\right\}.
\end{align*}
Then $\Phi(1, 1) = \left\{(0, 0, 0)^\top\right\}$, $\Phi(1, 2) = \Phi(1, 1) \cup \left\{(1, 0, 0)^\top\right\}$, 

$\Phi(2, 2) = \Phi(1, 2) \cup \left\{(0, 1, 0)^\top, (0, 0, 1)^\top\right\} = \Phi(2, 3)$,

\noindent cf. Figure~\ref{Fig_subbadd}, and we have $\varphi(1, 2) = -3 < -2 = \varphi(2, 2)$ as well as

$\varphi(1, 1) + \varphi(1, 2) = 0 - 3 = -3 < -2 = \varphi(2, 3)$.
\begin{figure}[ht]
\centering
\vspace{-0.15 in}
\begin{tikzpicture}[scale=0.7]
\filldraw (-0.5, -0.5) -- (0.33, -0.33) -- (1, 0) -- (0, 1);
\draw[dotted] (2, 0) -- (0, 1);
\draw[dotted] (2, 0) -- (-0.5, -0.5);
\path [pattern={Dots}] (2, 0) -- (0, 1) -- (1, 0) -- (0.33, -0.33) -- (2, 0);
\draw[dotted] (1, 0) -- (-1, -1);
\draw[dotted] (0, 1) -- (-1, -1);
\path [pattern={Dots}] (-1, -1) -- (0, 1) -- (-0.5, -0.5) -- (0.33, -0.33) -- (-1, -1);
\draw[thick, dashed] (2.5, 0.75) -- (0, 2) -- (-1, -1) -- (2.1, -0.38) -- (2.5, 0.75);
\path [pattern={Lines[angle=22.5,distance={8pt}]}] (2.5, 0.75) -- (0, 2) -- (-1, -1) -- (2.1, -0.38) -- (2.5, 0.75);
\path [pattern={Lines[angle=22.5,distance={8pt}]},pattern color=white] (-0.5, -0.5) -- (0.33, -0.33) -- (1, 0) -- (0, 1);
\draw[thick, ->] (-0.5, 0.5) -- (-0.8, 0.8);
\draw[thick, ->] (2, 1) -- (2.3, 1.15);
\draw[thick, ->] (1.5, -0.5) -- (1.85, -0.6);
\draw[thin, white] (0, 0) -- (1, 0);
\draw[thin, white] (0, 0) -- (0, 1);
\draw[thin, white] (0, 0) -- (-0.5, -0.5);
\draw[thick, ->] (1, 0) -- (2.5, 0) node[below] {$y_1$};
\draw[thick, ->] (0, 1) -- (0, 3.5) node[left] {$y_2$};
\draw[thick, ->] (-0.5, -0.5) -- (-1.5, -1.5) node[left] {$y_3$};
\filldraw[white] (1.85, -0.55) -- (2.15, -0.55) -- (2.15, -0.2) -- (1.85, -0.2);
\draw (2, 2pt) -- (2, -2pt) node[below] {$1$};
\draw (2pt, 2) -- (-2pt, 2) node[left] {$1$};
\draw (-1, -0.98) -- (-1, -1.02) node[below] {$1$};
\filldraw[white](0, 0)circle[radius=2pt];
\draw[thin] (0, 0)circle[radius=2pt];
\draw[dotted, white] (0, 1) -- (0.33, -0.33);
\path [pattern={Lines[angle=22.5,distance={4pt}]}] (-2, -2.55) -- (-1.6, -2.55) -- (-1.6, -2.15) -- (-2, -2.15);
\draw[thin] (-2, -2.55) -- (-1.6, -2.55) -- (-1.6, -2.15) -- (-2, -2.15) -- (-2, -2.55);
\draw[thin] (-1.8, -2.9)circle[radius=2pt];
\path [pattern={Lines[angle=22.5,distance={4pt}]}] (-2, -3.55) -- (-1.6, -3.55) -- (-1.6, -3.15) -- (-2, -3.15);
\draw[thin] (-2, -3.55) -- (-1.6, -3.55) -- (-1.6, -3.15) -- (-2, -3.15) -- (-2, -3.55);
\node[right] at (-1.6, -2.4) {$Wy \leq t$};
\node[right] at (-1.6, -2.9) {$\Phi(t)$};
\node[right] at (-1.6, -3.4) {$d^\top y \to \min!$};
\end{tikzpicture}
\begin{tikzpicture}[scale=0.7]
\filldraw (-1, -1) -- (2.33, -0.33) -- (3, 0) -- (2, 1) -- (0, 2);
\draw[dotted] (4, 0) -- (0, 2);
\draw[dotted] (4, 0) -- (-1, -1);
\path [pattern={Dots}] (4, 0) -- (2, 1) -- (3, 0) -- (2.33, -0.33) -- (4, 0);
\draw[dotted] (3, 0) -- (0, 3);
\draw[dotted] (3, 0) -- (-3, -3);
\draw[dotted] (0, 3) -- (-3, -3);
\path [pattern={Dots}] (-3, -3) -- (0, 3) -- (2, 1) -- (0, 2) -- (-1, -1) -- (2.33, -0.33) -- (-3, -3);
\path [pattern={Dots}] (-3, -3) -- (0, 3) -- (1, 2) -- (1.67, -0.67) -- (-3, -3);
\draw[thick, dashed] (4.5, 0.25) -- (0, 2.5) -- (-1.25, -1.25) -- (4.38, -0.13) -- (4.5, 0.25);
\path [pattern={Lines[angle=22.5,distance={8pt}]}] (4.5, 0.25) -- (0, 2.5) -- (-1.25, -1.25) -- (4.38, -0.13) -- (4.5, 0.25);
\path [pattern={Lines[angle=22.5,distance={8pt}]}, pattern color=white] (-1, -1) -- (2.33, -0.33) -- (3, 0) -- (2, 1) -- (0, 2);
\draw[thick, ->] (-0.63, 0.63) -- (-0.83, 0.83);
\draw[thick, ->] (2.5, 1.25) -- (2.8, 1.35);
\draw[thick, ->] (1.88, -0.62) -- (2.23, -0.72);
\draw[thin, white] (0, 0) -- (3, 0);
\draw[thin, white] (0, 0) -- (0, 2);
\draw[thin, white] (0, 0) -- (-1, -1);
\draw[thick, ->] (3, 0) -- (4.5, 0) node[below] {$y_1$};
\draw[thick, ->] (0, 2) -- (0, 3.5) node[left] {$y_2$};
\draw[thick, ->] (-1, -1) -- (-3.5, -3.5) node[left] {$y_3$};
\filldraw[white] (-0.45, 1.8) -- (-0.45, 2.2) -- (-0.2, 2.2) -- (-0.2, 1.8);
\filldraw[white] (-1.1, -1.55) -- (-1.1, -1.15) -- (-0.85, -1.15) -- (-0.85, -1.55);
\draw (2pt, 2) -- (-2pt, 2) node[left] {$1$};
\draw (-1, -0.98) -- (-1, -1.02) node[below] {$1$};
\filldraw[white] (0, 0)circle[radius=2pt]
    (2, 0)circle[radius=2pt]
    (0, 2)circle[radius=2pt]
    (-1, -1)circle[radius=2pt];
\draw[thin] (0, 0)circle[radius=2pt]
    (2, 0)circle[radius=2pt]
    (0, 2)circle[radius=2pt]
    (-1, -1)circle[radius=2pt];
\draw[dotted, white] (2, 1) -- (2.33, -0.33);
\end{tikzpicture}
\caption{Graph of the feasible points with a) $t = (1, 1)^\top$ and b) $t = (2, 3)^\top$ of Example~\ref{Ex_subadd}}
\label{Fig_subbadd}
\vspace{-0.15 in}
\end{figure}
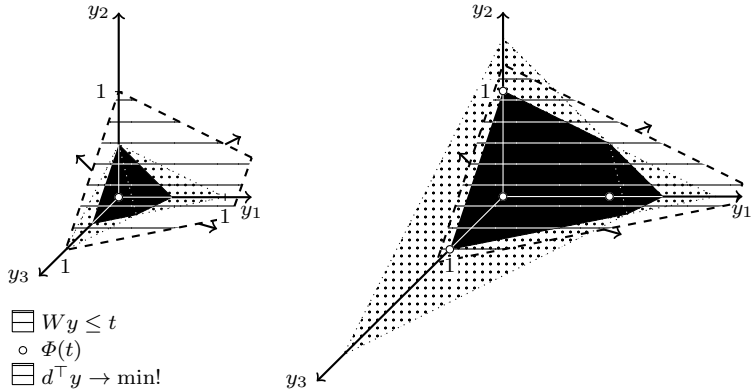
\end{example}

Based on \cite[Lemma 2]{ZhOe21} including the proof, here is a property of the function $\varphi$ that we will need later:
\begin{proposition} \label{Prop_phi}
It is $\varphi(t') \geq \varphi(t)$ for $t', t \in \mathrm{dom}\; \Phi$ such that $t' \leq t$ and $\psi(t') = \psi(t)$.
\end{proposition}

We shall work with the following assumption:
\begin{enumerate}[label=$\mathrm{(A\arabic*)}$,leftmargin=1cm]
    \item The set $Y \cap \mathbb{Z}^m$ is finite. \label{A1}
\end{enumerate}

\begin{lemma} \label{LemmaDomPhi}
We have $\mathrm{dom} \; \Phi = \bigcup_{y \in Y \cap \mathbb{Z}^m} \lbrace Wy \rbrace \oplus \mathbb{R}^s_{\geq 0}$, where $\oplus$ denotes the Minkowski sum, i.e. the domain of $\Phi$ admits a representation as a countable union of translated copies of the nonnegative orthant. In particular, $\mathrm{dom} \; \Phi$ is polyhedral, i.e. a finite union of polyhedra, cf. \cite[Def. 3.1]{De02}, whenever \ref{A1} holds true.
\end{lemma}

{\it Proof} This is an immediate consequence of the representation
\begin{align*}
    \mathrm{dom} \; \Phi 
    = \bigcup_{y \in Y \cap \mathbb{Z}^m} \lbrace t \; | \; Wy \leq t \rbrace = \bigcup_{y \in Y \cap \mathbb{Z}^m} \lbrace Wy \rbrace \oplus \mathbb{R}^s_{\geq 0}.\tag*{$\qed$}
\end{align*}

If $Y$ is a polyhedron and all integrality constraints are deleted from the follower's problem, it can be shown that the domain of $\varphi$ is empty or coincides with the set of parameters for which the lower level problem is feasible, i.e. $\mathrm{dom} \; \varphi = \mathrm{dom} \; \Phi$ (cf. \cite[Lemma 2.1]{BuClDe20}). The following example illustrates that this does not hold in the presence of integer variables:
\begin{example} \label{ExDomPhi}
Consider the case where $Y 
    := \lbrace (0,0)^\top \rbrace \cup \lbrace (1, k)^\top \; | \; k \in \mathbb{R} \rbrace$,
\noindent $m = 2$, $s = 1$, $W = [1\ 0]$ and $d = q = (0, 1)^\top$. It is easy to check that \mbox{$\mathrm{dom} \; \varphi = [0,1)$}, even though the lower level problem is feasible for any nonnegative right-hand side, cf. Figure~\ref{Fig_(1,k)}.
\end{example}

Example~\ref{ExDomPhi} also shows that the domain of $\varphi$ is not closed in general.
\begin{example} \label{Ex_Sqrt2}
Let $Y = \left\{y \in \mathbb{R}^3 \; | \; -\sqrt{2}y_1 + y_2 \leq 0, \, y_1 \geq 1, \, y_2 \geq 0, \, y_3 = 0\right\} \cup \left\{(0, 0, 1)^\top\right\}$, $W = [0\ 0\ 1]$, and $d = (\sqrt{2}, -1, -1)^\top$, then
\[
    \Psi(t) = 
    \begin{cases}
        \emptyset, & \text{if}\ t < 1\\
        \{(0, 0, 1)^\top \}, & \text{if}\ t \geq 1
    \end{cases}
\]
and the lower level infimal value $0$ is not attained for any $t \in [0, 1)$. In addition, we have
\[
    \Phi(t) = 
    \begin{cases}
        \emptyset, & \text{if}\ t < 0\\
        \left(Y \cap \mathbb{Z}^3\right) \setminus \{(0, 0, 1)^\top\}, & \text{if}\ t \in [0, 1)\\
        \{(0, 0, 1)^\top \}, & \text{if}\ t \geq 1,
    \end{cases}
\]
cf. Figure~\ref{Fig_sqrt2}, such that $\mathrm{dom} \; \Phi = [0, \infty) \supseteq [1, \infty) = \mathrm{dom} \, \Psi$.
\begin{figure}[ht]
\vspace{-0.15 in}
\begin{minipage}[t]{0.51\textwidth}
\centering
\begin{tikzpicture}[scale=0.7]
\path [pattern={Lines[angle=22.5,distance={8pt}]}] (-0.5, -2) -- (0.7, -2) -- (0.7, 1.5) -- (-0.5, 1.5);
\draw (0.7, -2.25) -- (0.7, 1.5);
\foreach \x in {1}
  \draw (\x,2pt) -- (\x,-2pt) node[below] {$\x$};
\foreach \y in {-1,1}
  \draw (2pt,\y) -- (-2pt,\y) node[left,fill=white] {$\y$};
\draw[thick, ->] (-0.5, 0) -- (1.5, 0) node[below] {$y_1$};
\draw[thick, ->] (0, -2) -- (0, 1.5) node[left] {$y_2$};
\draw[thick, ->] (-0.5, -2.1) -- (1.5, -2.1) node[below] {$t$};
\filldraw (0, -2.15)circle[radius=2pt];
    \draw[thick] (0, -2.15) -- (0.7, -2.15);
    \node[below] at (0.35, -2.15) {$t \in \mathrm{dom}\; \varphi$};
\filldraw[white] (0, 0)circle[radius=2pt]
(1, 0)circle[radius=2pt]
(1, 1)circle[radius=2pt]
(1, -1)circle[radius=2pt];
\draw[thin] (0, 0)circle[radius=2pt]
(1, 0)circle[radius=2pt]
(1, 1)circle[radius=2pt]
(1, -1)circle[radius=2pt];
\draw[thick] (-0.5, 0.5) -- (1.5, 0.5);
\draw[thick, ->] (0.5, 0.5) -- (0.5, 0.2);
\end{tikzpicture}
\hspace*{-1em}
\begin{tikzpicture}[scale=0.7]
\path [pattern={Lines[angle=22.5,distance={8pt}]}] (-0.5, -2) -- (1.7, -2) -- (1.7, 1.5) -- (-0.5, 1.5);
\draw (1.7, -2.25) -- (1.7, 1.5);
\foreach \x in {1}
  \draw (\x,2pt) -- (\x,-2pt) node[below,fill=white] {$\x$};
\foreach \y in {-1,1}
  \draw (2pt,\y) -- (-2pt,\y) node[left,fill=white] {$\y$};
\draw[thick, ->] (-0.5, 0) -- (2, 0) node[below] {$y_1$};
\draw[thick, ->] (0, -2) -- (0, 1.5) node[left] {$y_2$};
\draw[thick, ->] (-0.5, -2.1) -- (2, -2.1) node[below] {$t$};
\filldraw (0, -2.15)circle[radius=2pt];
    \draw[{-)},thick] (0, -2.15) -- (1, -2.15);
    \node[below] at (0.5, -2.15) {$t \in \mathrm{dom}\; \varphi$};
\draw[thick] (-0.5, -1.5) -- (2, -1.5);
\draw[thick, ->] (0.5, -1.5) -- (0.5, -1.8);
\filldraw[white] (0, 0)circle[radius=2pt]
(1, 0)circle[radius=2pt]
(1, 1)circle[radius=2pt]
(1, -1)circle[radius=2pt];
\draw[thin] (0, 0)circle[radius=2pt]
(1, 0)circle[radius=2pt]
(1, 1)circle[radius=2pt]
(1, -1)circle[radius=2pt];
\draw[thin] (2.3, 1.25)circle[radius=2pt];
\path [pattern={Lines[angle=22.5,distance={4pt}]}] (2.1, 0.5) -- (2.5, 0.5) -- (2.5, 1) -- (2.1, 1);
\draw[thin] (2.1, 0.5) -- (2.5, 0.5) -- (2.5, 1) -- (2.1, 1) -- (2.1, 0.5);
\draw[thin] (2.1, 0.4) -- (2.4, 0.1);
\draw[thin,->] (2.25, 0.25) -- (2.4, 0.4);
\node[right] at (2.5, 1.25) {$Y \cap \mathbb{Z}^2$};
\node[right] at (2.5, 0.75) {$Wy \leq t$};
\node[right] at (2.5, 0.25) {$d^\top y \!\to\! \min!$,};
\node[right] at (2.5, -.25) {$q^\top y \!\to\! \min!$};
\end{tikzpicture}
\caption{Graph of the feasible points and \mbox{$t \in \mathrm{dom}\; \varphi$} with a) $t \in (0, 1)$ and b) $t > 1$ of Example~\ref{ExDomPhi}}
\label{Fig_(1,k)}
\end{minipage}
\begin{minipage}[t]{0.48\textwidth}
\centering
\begin{tikzpicture}[scale=0.7]
\draw (1, 0) -- (1, 1.42) -- (2.47, 3.5);
\draw (1, 0.04) -- (4.5, 0.04);
\path [pattern={Lines[angle=22.5,distance={8pt}]}] (1, 0.04) -- (4.5, 0.04) -- (4.5, 3.5) -- (2.47, 3.5) -- (1, 1.42);
\draw[thick, ->] (0, 0) -- (4.5, 0) node[below] {$y_1$};
\draw[thick, ->] (0, 0) -- (0, 3.5) node[left] {$y_2$};
\foreach \x in {1,...,4}
    \draw (\x,2pt) -- (\x,-2pt) node[below] {$\x$};
\foreach \y in {1,...,3}
    \draw (2pt,\y) -- (-2pt,\y) node[left] {$\y$};
\draw[thick] (1, 0) -- (3.47, 3.5);
\draw[thick, ->] (2, 1.45) -- (1.7, 1.65);
\filldraw[white] (1, 0)circle[radius=2pt]
    (1, 1)circle[radius=2pt]
    (2, 0)circle[radius=2pt]
    (2, 1)circle[radius=2pt]
    (2, 2)circle[radius=2pt]
    (3, 0)circle[radius=2pt]
    (3, 1)circle[radius=2pt]
    (3, 2)circle[radius=2pt]
    (3, 3)circle[radius=2pt]
    (4, 0)circle[radius=2pt]
    (4, 1)circle[radius=2pt]
    (4, 2)circle[radius=2pt]
    (4, 3)circle[radius=2pt];
\draw[thin] (1, 0)circle[radius=2pt]
    (1, 1)circle[radius=2pt]
    (2, 0)circle[radius=2pt]
    (2, 1)circle[radius=2pt]
    (2, 2)circle[radius=2pt]
    (3, 0)circle[radius=2pt]
    (3, 1)circle[radius=2pt]
    (3, 2)circle[radius=2pt]
    (3, 3)circle[radius=2pt]
    (4, 0)circle[radius=2pt]
    (4, 1)circle[radius=2pt]
    (4, 2)circle[radius=2pt]
    (4, 3)circle[radius=2pt];
\path [pattern={Lines[angle=22.5,distance={4pt}]}] (4.6, 3) -- (5, 3) -- (5, 3.5) -- (4.6, 3.5);
\draw[thin] (4.6, 3) -- (5, 3) -- (5, 3.5) -- (4.6, 3.5) -- (4.6, 3);
\draw[thin] (4.8, 2.75)circle[radius=2pt];
\draw[thin] (4.6, 2.4) -- (4.9, 2.1);
\draw[thin,->] (4.75, 2.25) -- (4.9, 2.4);
\node[right] at (5, 3.25) {$Y \subseteq \mathbb{R}^2$};
\node[right] at (5, 2.75) {$Y \cap \mathbb{Z}^2$};
\node[right] at (5, 2.25) {$d^\top y \!\to\! \min!$};
\end{tikzpicture}
\caption{Graph of $\Phi(t)$ for \mbox{$t \in [0, 1)$} with \mbox{$y_3 = 0$} of Example~\ref{Ex_Sqrt2}}
\label{Fig_sqrt2}
\end{minipage}
\vspace{-0.15 in}
\end{figure}

If we instead set $d = (0, 0, 0)^\top$ and $q = (\sqrt{2}, -1, -1)^\top$, it is easy to see that the infimal value of the optimization problem defining $\varphi(t)$ is finite but not attained for any $t \in [0,1)$. Here we have $\mathrm{dom} \; \Phi = [0, \infty) = \mathrm{dom} \, \Psi$, but $\mathrm{dom} \; \psi = [0, \infty) \supseteq [1, \infty) = \mathrm{dom} \, \varphi$.
\end{example}

In view of Lemma~\ref{LemmaDomPhi} it is desirable to identify situations in which the domains of $\varphi$ and $\Phi$ coincide:
\begin{lemma} \label{LemmaA1DomainsPolyhedral}
Assume \ref{A1}, then the domains of $\varphi, \psi, \Psi$ and $\Phi$ coincide and are polyhedral.
\end{lemma}

{\it Proof} The lower level problem as well as the optimization problem defining $\varphi$ have at most $|Y \cap \mathbb{Z}^m| < \infty$ feasible points, which implies the problems are solvable whenever they are feasible. Hence \mbox{$\mathrm{dom} \; \varphi = \mathrm{dom} \; \Psi = \mathrm{dom} \; \psi = \mathrm{dom} \; \Phi$} and the statement follows directly from Lemma~\ref{LemmaDomPhi}.\qed

\smallskip
For the subsequent analysis, let us assume that $Y \cap \mathbb{Z}^m$ is nonempty and \ref{A1} holds, i.e. $Y \cap \mathbb{Z}^m 
    = \left\{\bar{y}^1,\ldots,\bar{y}^N\right\}$ for some $N \in \mathbb{N}$. The following example illustrates that the domain of $\varphi$ may be nonconvex if $s \geq 2$. 
\begin{example} \label{Ex_2}
Figure~\ref{Fig_Y} depicts a bounded set containing $N = 6$ integral points.
\begin{enumerate}[leftmargin=2em, label={\alph*)}]
    \item For $W = [1\ 1]$, we have $s = 1$ and $\mathrm{dom}\;\varphi$ is shown in Figure~\ref{Fig_domphi_a}.
\begin{figure}[ht]
\vspace{-0.15 in}
\begin{minipage}[t]{0.33\textwidth}
\centering
\begin{tikzpicture}[scale=0.7]
\path[pattern={Lines[angle=22.5,distance={8pt}]}] (3, 1) to [thin, closed, curve through =
    {(4.2, 2) (3.2, 3.2) 
    (1.9, 3.1) (2, 2.5) 
    (2.4, 1.6) (1.5, 1.8)
    (0.9, 0.9) (2, 1.2)}] (3, 1);
\draw[thin] (3, 1) to [thin, closed, curve through =
    {(4.2, 2) (3.2, 3.2) 
    (1.9, 3.1) (2, 2.5) 
    (2.4, 1.6) (1.5, 1.8)
    (0.9, 0.9) (2, 1.2)}] (3, 1);
\node [rectangle,minimum width=0.4cm,minimum height=0.4cm,fill=white] at (2, 2.65) {};
\node [rectangle,minimum width=0.4cm,minimum height=0.4cm,fill=white] at (3, 2.65) {};
\node [rectangle,minimum width=0.4cm,minimum height=0.4cm,fill=white] at (3, 1.65) {};
\node [rectangle,minimum width=0.4cm,minimum height=0.4cm,fill=white] at (4, 1.65) {};
\node [rectangle,minimum width=0.4cm,minimum height=0.4cm,fill=white] at (1, 0.65) {};
\filldraw[white] (1, 1)circle[radius=2pt]
    (3, 1)circle[radius=2pt]
    (3, 2)circle[radius=2pt]
    (4, 2)circle[radius=2pt]
    (2, 3)circle[radius=2pt]
    (3, 3)circle[radius=2pt];
\draw[thin] (1, 1)circle[radius=2pt] node[below] {$\bar{y}^5$}
    (3, 1)circle[radius=2pt] node[below] {$\bar{y}^6$}
    (3, 2)circle[radius=2pt] node[below] {$\bar{y}^3$}
    (4, 2)circle[radius=2pt] node[below] {$\bar{y}^4$}
    (2, 3)circle[radius=2pt] node[below] {$\bar{y}^1$}
    (3, 3)circle[radius=2pt] node[below] {$\bar{y}^2$};
\draw[thick, ->] (0, 0) -- (4.5, 0) node[below] {$y_1$};
\draw[thick, ->] (0, 0) -- (0, 3.5) node[left] {$y_2$};
\foreach \x in {1,...,4}
    \draw (\x,2pt) -- (\x,-2pt) node[below] {$\x$};
\foreach \y in {1,...,3}
    \draw (2pt,\y) -- (-2pt,\y) node[left] {$\y$};
    \node[white,below] at (4.5, -0.4) {$1$};
\end{tikzpicture}
\caption{The set of integral points in a bounded subset $Y$ of $\mathbb{R}^2$ in Example~\ref{Ex_2}}
\label{Fig_Y}
\end{minipage}
\begin{minipage}[t]{0.66\textwidth}
\centering
\begin{tikzpicture}[scale=0.7]
\draw[thick, ->] (0, 0) -- (6.5, 0);
\node[below] at (6.5, -0.4) {$t$};
\foreach \x in {1,...,6}
    \draw (\x,2pt) -- (\x,-2pt);
    \node[below] at (1, -0.05) {$1$};
    \node[below] at (2, -0.1) {$2$};
    \node[below] at (3, -0.1) {$3$};
    \node[below] at (4, -0.2) {$4$};
    \node[below] at (5, -0.3) {$5$};
    \node[below] at (6, -0.4) {$6$};
\draw[dotted, thin] (1, 1) -- (2, 0);
    \filldraw (2, -0.1)circle[radius=2pt];
    \draw[thick] (2, -0.1) -- (7, -0.1);
\draw[dashed, thin] (3, 1) -- (4, 0);
    \filldraw (4, -0.2)circle[radius=2pt];
    \draw[thick] (4, -0.2) -- (7, -0.2);
\draw[dotted, thin] (2, 3) -- (5, 0);
    \filldraw (5, -0.3)circle[radius=2pt];
    \draw[thick] (5, -0.3) -- (7, -0.3);
\draw[dashed, thin] (3, 3) -- (6, 0);
    \filldraw (6, -0.4)circle[radius=2pt];
    \draw[thick] (6, -0.4) -- (7, -0.4);
\filldraw[white] (1, 1)circle[radius=2pt]
    (3, 1)circle[radius=2pt]
    (3, 2)circle[radius=2pt]
    (4, 2)circle[radius=2pt]
    (2, 3)circle[radius=2pt]
    (3, 3)circle[radius=2pt];
\draw[thin] (1, 1)circle[radius=2pt] node[below] {$\bar{y}^5$}
    (3, 1)circle[radius=2pt] node[below] {$\bar{y}^6$}
    (3, 2)circle[radius=2pt] node[below] {$\bar{y}^3$}
    (4, 2)circle[radius=2pt] node[below] {$\bar{y}^4$}
    (2, 3)circle[radius=2pt] node[below] {$\bar{y}^1$}
    (3, 3)circle[radius=2pt] node[below] {$\bar{y}^2$};
\draw[thin] (7.3, 3.25)circle[radius=2pt] node {};
\draw[dotted, thin] (7.5, 2.75) -- (7, -0.1);
\draw[dashed, thin] (7.5, 1.75) -- (7, -0.2);
\draw[dotted, thin] (7.5, 1.25) -- (7, -0.3);
\draw[dashed, thin] (7.5, 0.75) -- (7, -0.4);
\node[right] at (7.5, 3.25) {$Y \cap \mathbb{Z}^2$};
\node[right] at (7.5, 2.75) {$\{t | \bar{y}^i_1 \!+\! \bar{y}^i_2 \!=\! 2 \!\leq\! t\}$};
\node[right] at (7.55, 2.25) {$= \mathrm{dom}\;\varphi$};
\node[right] at (7.5, 1.75) {$\{t | \bar{y}^i_1 \!+\! \bar{y}^i_2 \!=\! 4 \!\leq\! t\}$};
\node[right] at (7.5, 1.25) {$\{t | \bar{y}^i_1 \!+\! \bar{y}^i_2 \!=\! 5 \!\leq\! t\}$};
\node[right] at (7.5, 0.75) {$\{t | \bar{y}^i_1 \!+\! \bar{y}^i_2 \!=\! 6 \!\leq\! t\}$};
\end{tikzpicture}
\caption{The domain of $\varphi$ in Example~\ref{Ex_2}~a)}
\label{Fig_domphi_a}
\end{minipage}
\vspace{-0.15 in}
\end{figure}
    \item As in Example~\ref{Ex_1}, we use $W = \bigl[\begin{smallmatrix}
        -1 & 0\\
        0 & 1
    \end{smallmatrix}\bigr]$ and have $s = 2$. Then $\mathrm{dom}\;\varphi$ is shown in Figure~\ref{Fig_domphi_b}. 
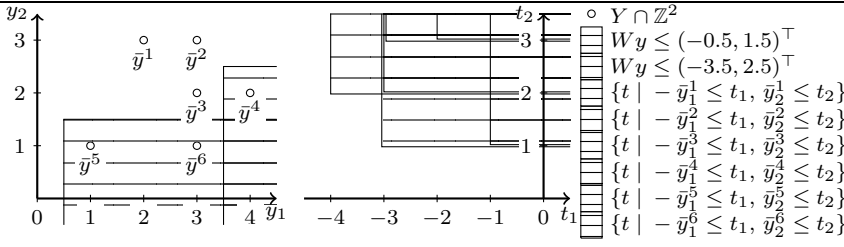
\begin{figure}[ht]
\centering
\vspace{-0.15 in}
\begin{tikzpicture}[scale=0.7]
\path [pattern={Lines[angle=67.5,distance={8pt}]}] (4.5, 1.5) -- (0.5, 1.5) -- (0.5, -0.5) -- (4.5, -0.5);
\draw (4.5, 1.5) -- (0.5, 1.5) -- (0.5, -0.5);
\path [pattern={Lines[angle=-67.5,distance={8pt}]}] (4.5, 2.5) -- (3.5, 2.5) -- (3.5, -0.5) -- (4.5, -0.5);
\draw (4.5, 2.5) -- (3.5, 2.5) -- (3.5, -0.5);
\node [rectangle,minimum width=0.4cm,minimum height=0.4cm,fill=white] at (3, 1.65) {};
\node [rectangle,minimum width=0.4cm,minimum height=0.4cm,fill=white] at (4, 1.65) {};
\node [rectangle,minimum width=0.4cm,minimum height=0.4cm,fill=white] at (1, 0.65) {};
\node [rectangle,minimum width=0.4cm,minimum height=0.4cm,fill=white] at (3, 0.65) {};
\node [rectangle,minimum width=0.3cm,minimum height=0.25cm,fill=white] at (4.45, -.3) {};
\draw[thick, ->] (0, 0) -- (4.5, 0) node[below] {$y_1$};
\draw[thick, ->] (0, 0) -- (0, 3.5) node[left] at (0, 3.5) {$y_2$};
\foreach \x in {0,...,4}
    \draw (\x,2pt) -- (\x,-2pt) node[below, fill=white] {$\x$};
\foreach \y in {1,...,3}
    \draw (2pt,\y) -- (-2pt,\y) node[left] {$\y$};
\filldraw[white] (1, 1)circle[radius=2pt]
    (3, 1)circle[radius=2pt]
    (3, 2)circle[radius=2pt]
    (4, 2)circle[radius=2pt]
    (2, 3)circle[radius=2pt]
    (3, 3)circle[radius=2pt];
\draw[thin] (1, 1)circle[radius=2pt] node[below] {$\bar{y}^5$}
    (3, 1)circle[radius=2pt] node[below] {$\bar{y}^6$}
    (3, 2)circle[radius=2pt] node[below] {$\bar{y}^3$}
    (4, 2)circle[radius=2pt] node[below] {$\bar{y}^4$}
    (2, 3)circle[radius=2pt] node[below] {$\bar{y}^1$}
    (3, 3)circle[radius=2pt] node[below] {$\bar{y}^2$};
    \node[white,below] at (4, -0.32) {$1$};
\end{tikzpicture}
\begin{tikzpicture}[scale=0.7]
\path [pattern={Lines[angle=-45,distance={8pt}]}] (-4, 2) -- (0.5, 2) -- (0.5, 3.5) -- (-4, 3.5);
\path [pattern={Lines[angle=0,distance={8pt}]}] (-3, 1) -- (0.5, 1) -- (0.5, 3.5) -- (-3, 3.5);
\path [pattern={Lines[angle=-22.5,distance={8pt},yshift=-3pt]}] (-1, 1) -- (0.5, 1) -- (0.5, 3.5) -- (-1, 3.5);
\path [pattern={Lines[angle=90,distance={8pt},yshift=-3pt]}] (-3, 2) -- (0.5, 2) -- (0.5, 3.5) -- (-3, 3.5);
\path [pattern={Lines[angle=45,distance={8pt}]}] (-3, 3) -- (0.5, 3) -- (0.5, 3.5) -- (-3, 3.5);
\path [pattern={Lines[angle=22.5,distance={8pt},yshift=-3pt]}] (-2, 3) -- (0.5, 3) -- (0.5, 3.5) -- (-2, 3.5);
\draw (-4, 3.5) -- (-4, 1.98) -- (0.5, 1.98);
\draw (-3.04, 3.5) -- (-3.04, 0.98) -- (0.5, 0.98);
\draw (-1, 3.5) -- (-1, 1.02) -- (0.5, 1.02);
\draw (-3, 3.5) -- (-3, 2.02) -- (0.5, 2.02);
\draw (-2.96, 3.5) -- (-2.96, 2.98) -- (0.5, 2.98);
\draw (-2, 3.5) -- (-2, 3.02) -- (0.5, 3.02);
\node [rectangle,minimum width=0.25cm,minimum height=0.4cm,fill=white] at (-0.28, 1) {};
\node [rectangle,minimum width=0.28cm,minimum height=0.4cm,fill=white] at (-0.26, 2) {};
\node [rectangle,minimum width=0.28cm,minimum height=0.4cm,fill=white] at (-0.26, 3) {};
\node [rectangle,minimum width=0.28cm,minimum height=0.3cm,fill=white] at (-0.3, 3.45) {};
\draw[thick, ->] (-4.5, 0) -- (0.5, 0) node[below] {$t_1$};
\draw[thick, ->] (0, 0) -- (0, 3.5) node[left] {$t_2$};
\foreach \x in {-4,...,0}
    \draw (\x,2pt) -- (\x,-2pt) node[below] {$\x$};
\foreach \y in {1,...,3}
    \draw (2pt,\y) -- (-2pt,\y) node[left] {$\y$};
\draw[thin] (0.9, 3.5)circle[radius=2pt] node[below] {};
\path [pattern={Lines[angle=67.5,distance={4pt}]}] (0.7, 2.75) -- (1.1, 2.75) -- (1.1, 3.25) -- (0.7, 3.25);
\draw[thin] (0.7, 2.75) -- (1.1, 2.75) -- (1.1, 3.25) -- (0.7, 3.25) -- (0.7, 2.75);
\path [pattern={Lines[angle=-67.5,distance={4pt}]}] (0.7, 2.25) -- (1.1, 2.25) -- (1.1, 2.75) -- (0.7, 2.75);
\draw[thin] (0.7, 2.25) -- (1.1, 2.25) -- (1.1, 2.75) -- (0.7, 2.75) -- (0.7, 2.25);
\path [pattern={Lines[angle=22.5,distance={4pt}]}] (0.7, 1.75) -- (1.1, 1.75) -- (1.1, 2.25) -- (0.7, 2.25);
\draw[thin] (0.7, 1.75) -- (1.1, 1.75) -- (1.1, 2.25) -- (0.7, 2.25) -- (0.7, 1.75);
\path [pattern={Lines[angle=45,distance={4pt}]}] (0.7, 1.25) -- (1.1, 1.25) -- (1.1, 1.75) -- (0.7, 1.75);
\draw[thin] (0.7, 1.25) -- (1.1, 1.25) -- (1.1, 1.75) -- (0.7, 1.75) -- (0.7, 1.25);
\path [pattern={Lines[angle=90,distance={4pt}]}] (0.7, 0.75) -- (1.1, 0.75) -- (1.1, 1.25) -- (0.7, 1.25);
\draw[thin] (0.7, 0.75) -- (1.1, 0.75) -- (1.1, 1.25) -- (0.7, 1.25) -- (0.7, 0.75);
\path [pattern={Lines[angle=-45,distance={4pt}]}] (0.7, 0.25) -- (1.1, 0.25) -- (1.1, 0.75) -- (0.7, 0.75);
\draw[thin] (0.7, 0.25) -- (1.1, 0.25) -- (1.1, 0.75) -- (0.7, 0.75) -- (0.7, 0.25);
\path [pattern={Lines[angle=-22.5,distance={4pt}]}] (0.7, -.25) -- (1.1, -.25) -- (1.1, 0.25) -- (0.7, 0.25);
\draw[thin] (0.7, -.25) -- (1.1, -.25) -- (1.1, 0.25) -- (0.7, 0.25) -- (0.7, -.25);
\path [pattern={Lines[angle=0,distance={4pt}]}] (0.7, -.75) -- (1.1, -.75) -- (1.1, -.25) -- (0.7, -.25);
\draw[thin] (0.7, -.75) -- (1.1, -.75) -- (1.1, -.25) -- (0.7, -.25) -- (0.7, -0.75);
\node[right] at (1.1, 3.5) {$Y \cap \mathbb{Z}^2$};
\node[right] at (1.1, 3) {$Wy \leq (-0.5, 1.5)^\top$};
\node[right] at (1.1, 2.5) {$Wy \leq (-3.5, 2.5)^\top$};
\node[right] at (1.1, 2) {$\{t\; | \; -\bar{y}^1_1 \leq t_1, \, \bar{y}^1_2 \leq t_2\}$};
\node[right] at (1.1, 1.5) {$\{t\; | \; -\bar{y}^2_1 \leq t_1, \, \bar{y}^2_2 \leq t_2\}$};
\node[right] at (1.1, 1) {$\{t\; | \; -\bar{y}^3_1 \leq t_1, \, \bar{y}^3_2 \leq t_2\}$};
\node[right] at (1.1, 0.5) {$\{t\; | \; -\bar{y}^4_1 \leq t_1, \, \bar{y}^4_2 \leq t_2\}$};
\node[right] at (1.1, 0) {$\{t\; | \; -\bar{y}^5_1 \leq t_1, \, \bar{y}^5_2 \leq t_2\}$};
\node[right] at (1.1, -.5) {$\{t\; | \; -\bar{y}^6_1 \leq t_1, \, \bar{y}^6_2 \leq t_2\}$};
\end{tikzpicture}
\caption{Graph of a) the feasible points for $t = (-0.5, 1.5)$ and $t = (-3.5, 2.5)$, and b) of $\mathrm{dom}\, \varphi$ \footnotesize{(boundaries of the sets are plotted side by side or better visibility)} of Example~\ref{Ex_2}~b)}
\label{Fig_domphi_b}
\vspace{-0.15 in}
\end{figure}
\end{enumerate}
\end{example}

As preliminary work for the analysis in Section~\ref{Sec_Model}, let us enumerate the nonempty subsets of $Y \cap \mathbb{Z}^m$: $\{A \subseteq Y \cap \mathbb{Z}^m \; | \; A \neq \emptyset \} = \left\{Y^1, \ldots, Y^M\right\}$ with $M = 2^N-1$, which allows us to define
\allowdisplaybreaks
\begin{align*}
    V^k 
    := &\left\{t \; | \; \left\{y \in Y \cap \mathbb{Z}^m \; | \; Wy \leq t\right\} = Y^k\right\} \\
    = &\bigcap_{y \in Y^k} \{t \; | \; Wy \leq t\} \cap
    \bigcap_{\hat{y} \in \left(Y \cap \mathbb{Z}^m\right) \setminus Y^k} \bigcup_{j = 1}^s \left\{t \; | \; e_j^\top W\hat{y} > t_j\right\}
\end{align*}
for any $k \in \lbrace 1, \ldots, M \rbrace$. The first intersection above guarantees that all points from $Y^k$ are feasible, while the second one rules out all other elements in $Y \cap \mathbb{Z}^m$. As illustrated in the continuation of Example~\ref{Ex_2} below, the sets $V^k$ form a partion of the domain of $\varphi$. This is indeed always true:
\begin{lemma} \label{Lemma_partition}
Assume \ref{A1}, then the family $\mathcal{V} := \{V^k \; | \; V^k \neq \emptyset, \, k = 1, \ldots, M\}$ forms a partition of $\mathrm{dom}\; \varphi$.
\end{lemma}

{\it Proof} Lemma~\ref{LemmaA1DomainsPolyhedral} yields $\mathrm{dom} \; \varphi = \lbrace t \; | \; \exists k \in \lbrace 1, \dots, M \rbrace \; : \; t \in V^k \rbrace$ and the statement follows directly from the fact that the sets in $\mathcal{V}$ are the equivalence classes w.r.t. the equivalence relation where $t \sim t'$ iff
\begin{align*}
    \lbrace y \in Y \cap \mathbb{Z}^m \; | \; Wy \leq t \rbrace = \lbrace y \in Y \cap \mathbb{Z}^m \; | \; Wy \leq t' \rbrace.\tag*{$\qed$}
\end{align*}

By definition, the family $\mathcal{V}$ may have up to $2^{|Y \cap \mathbb{Z}^m|}$ elements. However, the following observation shows that the cardinality of $\mathcal{V}$ is usually much smaller:
\begin{corollary} \label{Cor_emptyset}
Assume \ref{A1} and let $k = 1, \ldots, M$ as well as $i, \iota \in \{1, \ldots, N\}$ be such that $\bar{y}^i \notin Y^k$, $\bar{y}^\iota \in Y^k$ and $W\bar{y}^i \leq W\bar{y}^\iota$, then $V^k = \emptyset$.
\end{corollary}

{\it Proof} Suppose that there is some $t \in V^k$, then $W\bar{y}^i \leq W\bar{y}^\iota \leq t$ implies $\bar{y}^i \notin Y^k$, which contradicts the assumptions. \qed

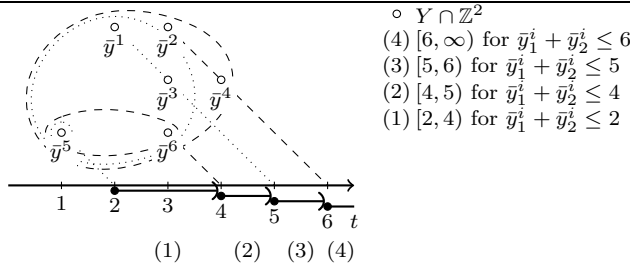
\begin{figure}[ht]
\centering
\vspace{-0.15 in}
\begin{tikzpicture}[scale=0.7]
\draw[thick, ->] (0, 0) -- (6.5, 0);
\node[below] at (6.5, -0.4) {$t$};
\foreach \x in {1,...,6}
    \draw (\x,2pt) -- (\x,-2pt);
    \node[below] at (1, -0.05) {$1$};
    \node[below] at (2, -0.1) {$2$};
    \node[below] at (3, -0.1) {$3$};
    \node[below] at (4, -0.2) {$4$};
    \node[below] at (5, -0.3) {$5$};
    \node[below] at (6, -0.4) {$6$};
\draw[dotted] (1, 1.2) to [closed, curve through =
    {(0.8, 1) (1, 0.8)}] (1.2, 1);
\draw[dashed] (1, 1.3) to [closed, curve through =
    {(0.7, 1) (1, 0.7) (3, 0.8) (3.2, 1) (3, 1.2)}] (1, 1.3);
\draw[dotted] (0.6, 1) to [closed, curve through =
    {(1, 0.6) (2.5, 0.5) (3, 0.7) (3.3, 1) (3, 2.8) (2, 3.2) (1, 2.8)}] (0.6, 1);
\draw[dashed] (0.5, 1) to [closed, curve through =
    {(1, 0.5) (3, 0.6) (4.2, 2.2) (3, 3.2) (2, 3.3)}] (0.5, 1);
\draw[dotted] (1, 1) -- (2, 0);
    \filldraw (2, -0.1)circle[radius=2pt];
    \draw[{-)}, thick] (2, -0.1) -- (3.95, -0.1);
\draw[dashed] (3, 1) -- (4, 0);
    \filldraw (4, -0.2)circle[radius=2pt];
    \draw[{-)}, thick] (4, -0.2) -- (4.95, -0.2);
\draw[dotted] (2, 3) -- (5, 0);
    \filldraw (5, -0.3)circle[radius=2pt];
    \draw[{-)}, thick] (5, -0.3) -- (5.95, -0.3);
\draw[dashed] (3, 3) -- (6, 0);
    \filldraw (6, -0.4)circle[radius=2pt];
    \draw[thick] (6, -0.4) -- (6.5, -0.4);
\node [rectangle,minimum width=0.4cm,minimum height=0.4cm,fill=white] at (2, 2.65) {};
\node [rectangle,minimum width=0.4cm,minimum height=0.4cm,fill=white] at (3, 2.65) {};
\node [rectangle,minimum width=0.4cm,minimum height=0.4cm,fill=white] at (3, 1.65) {};
\node [rectangle,minimum width=0.4cm,minimum height=0.4cm,fill=white] at (4, 1.65) {};
\node [rectangle,minimum width=0.4cm,minimum height=0.4cm,fill=white] at (1, 0.65) {};
\node [rectangle,minimum width=0.4cm,minimum height=0.4cm,fill=white] at (3, 0.65) {};
\filldraw[white] (1, 1)circle[radius=2pt]
    (3, 1)circle[radius=2pt]
    (3, 2)circle[radius=2pt]
    (4, 2)circle[radius=2pt]
    (2, 3)circle[radius=2pt]
    (3, 3)circle[radius=2pt];
\draw[thin] (1, 1)circle[radius=2pt] node[below] {$\bar{y}^5$}
    (3, 1)circle[radius=2pt] node[below] {$\bar{y}^6$}
    (3, 2)circle[radius=2pt] node[below] {$\bar{y}^3$}
    (4, 2)circle[radius=2pt] node[below] {$\bar{y}^4$}
    (2, 3)circle[radius=2pt] node[below] {$\bar{y}^1$}
    (3, 3)circle[radius=2pt] node[below] {$\bar{y}^2$};
\node at (6.25, -1.25) {$(4)$};
\node at (5.5, -1.25) {$(3)$};
\node at (4.5, -1.25) {$(2)$};
\node at (3, -1.25) {$(1)$};
\draw[thin] (7.3, 3.25)circle[radius=2pt];
\node at (7.3, 2.75) {$(4)$};
\node at (7.3, 2.25) {$(3)$};
\node at (7.3, 1.75) {$(2)$};
\node at (7.3, 1.25) {$(1)$};
\node[right] at (7.5, 3.25) {$Y \cap \mathbb{Z}^2$};
\node[right] at (7.5, 2.75) {$[6, \infty)$ for $\bar{y}^i_1 + \bar{y}^i_2 \leq 6$};
\node[right] at (7.5, 2.25) {$[5, 6)$ for $\bar{y}^i_1 + \bar{y}^i_2 \leq 5$};
\node[right] at (7.5, 1.75) {$[4, 5)$ for $\bar{y}^i_1 + \bar{y}^i_2 \leq 4$};
\node[right] at (7.5, 1.25) {$[2, 4)$ for $\bar{y}^i_1 + \bar{y}^i_2 \leq 2$};
\end{tikzpicture}
\caption{Graph of the sets $V^k$ of Example~\ref{Ex_2}~a)}
\label{Fig_Yk_a}
\vspace{-0.15 in}
\end{figure}

\begin{figure}[ht]
\centering
\vspace{-0.15 in}
\begin{tikzpicture}[scale=0.7]
\path [pattern={Lines[angle=22.5,distance={8pt}]}] (-4, 2) -- (-3, 2) -- (-3, 3.5) -- (-4, 3.5);
\path [pattern={Lines[angle=45,distance={8pt}]}] (-3, 1) -- (-1, 1) -- (-1, 2) -- (-3, 2);
\path [pattern={Lines[angle=67.5,distance={8pt}]}] (-3, 2) -- (-1, 2) -- (-1, 3) -- (-3, 3);
\path [pattern={Lines[angle=90,distance={8pt}]}] (-3, 3) -- (-2, 3) -- (-2, 3.5) -- (-3, 3.5);
\path [pattern={Lines[angle=-67.5,distance={8pt}]}] (-2, 3) -- (-1, 3) -- (-1, 3.5) -- (-2, 3.5);
\path [pattern={Lines[angle=-45,distance={8pt}]}] (-1, 1) -- (0.5, 1) -- (0.5, 2) -- (-1, 2);
\path [pattern={Lines[angle=-22.5,distance={8pt}]}] (-1, 2) -- (0.5, 2) -- (0.5, 3) -- (-1, 3);
\path [pattern={Lines[angle=0,distance={8pt}]}] (-1, 3) -- (0.5, 3) -- (0.5, 3.5) -- (-1, 3.5);
\draw (-4, 3.5) -- (-4, 2) -- (-3, 2) -- (-3, 3.5);
\draw (-3, 2) -- (-3, 1) -- (0.5, 1);
\draw (-3, 2) -- (0.5, 2);
\draw (-3, 3) -- (0.5, 3);
\draw (-2, 3) -- (-2, 3.5);
\draw (-1, 1) -- (-1, 3.5);
\node [rectangle,minimum width=0.25cm,minimum height=0.4cm,fill=white] at (-0.28, 1) {};
\node [rectangle,minimum width=0.28cm,minimum height=0.4cm,fill=white] at (-0.26, 2) {};
\node [rectangle,minimum width=0.28cm,minimum height=0.4cm,fill=white] at (-0.26, 3) {};
\node [rectangle,minimum width=0.28cm,minimum height=0.3cm,fill=white] at (-0.3, 3.45) {};
\node [rectangle,minimum width=0.4cm,minimum height=0.4cm,fill=white] at (2, 2.65) {};
\node [rectangle,minimum width=0.4cm,minimum height=0.4cm,fill=white] at (3, 2.65) {};
\node [rectangle,minimum width=0.4cm,minimum height=0.4cm,fill=white] at (3, 1.65) {};
\node [rectangle,minimum width=0.4cm,minimum height=0.4cm,fill=white] at (4, 1.65) {};
\node [rectangle,minimum width=0.4cm,minimum height=0.4cm,fill=white] at (1, 0.65) {};
\node [rectangle,minimum width=0.4cm,minimum height=0.4cm,fill=white] at (3, 0.65) {};
\draw[thick, ->] (-4.5, 0) -- (0.5, 0) node[below] {$t_1$};
\draw[thick, ->] (0, 0) -- (0, 3.5) node[left] {$t_2$};
\foreach \x in {-4,...,0}
    \draw (\x,2pt) -- (\x,-2pt) node[below] {$\x$};
\foreach \y in {1,...,3}
    \draw (2pt,\y) -- (-2pt,\y) node[left] {$\y$};
\draw[thin] (0.9, 3.5)circle[radius=2pt];
\path [pattern={Lines[angle=22.5,distance={4pt}]}] (0.7, 2.75) -- (1.1, 2.75) -- (1.1, 3.25) -- (0.7, 3.25);
\draw[thin] (0.7, 2.75) -- (1.1, 2.75) -- (1.1, 3.25) -- (0.7, 3.25) -- (0.7, 2.75);
\path [pattern={Lines[angle=45,distance={4pt}]}] (0.7, 2.25) -- (1.1, 2.25) -- (1.1, 2.75) -- (0.7, 2.75);
\draw[thin] (0.7, 2.25) -- (1.1, 2.25) -- (1.1, 2.75) -- (0.7, 2.75) -- (0.7, 2.25);
\path [pattern={Lines[angle=67.5,distance={4pt}]}] (0.7, 1.75) -- (1.1, 1.75) -- (1.1, 2.25) -- (0.7, 2.25);
\draw[thin] (0.7, 1.75) -- (1.1, 1.75) -- (1.1, 2.25) -- (0.7, 2.25) -- (0.7, 1.75);
\path [pattern={Lines[angle=90,distance={4pt}]}] (0.7, 1.25) -- (1.1, 1.25) -- (1.1, 1.75) -- (0.7, 1.75);
\draw[thin] (0.7, 1.25) -- (1.1, 1.25) -- (1.1, 1.75) -- (0.7, 1.75) -- (0.7, 1.25);
\path [pattern={Lines[angle=-67.5,distance={4pt}]}] (0.7, 0.75) -- (1.1, 0.75) -- (1.1, 1.25) -- (0.7, 1.25);
\draw[thin] (0.7, 0.75) -- (1.1, 0.75) -- (1.1, 1.25) -- (0.7, 1.25) -- (0.7, 0.75);
\path [pattern={Lines[angle=-45,distance={4pt}]}] (0.7, 0.25) -- (1.1, 0.25) -- (1.1, 0.75) -- (0.7, 0.75);
\draw[thin] (0.7, 0.25) -- (1.1, 0.25) -- (1.1, 0.75) -- (0.7, 0.75) -- (0.7, 0.25);
\path [pattern={Lines[angle=-22.5,distance={4pt}]}] (0.7, -0.25) -- (1.1, -0.25) -- (1.1, 0.25) -- (0.7, 0.25);
\draw[thin] (0.7, -0.25) -- (1.1, -0.25) -- (1.1, 0.25) -- (0.7, 0.25) -- (0.7, -0.25);
\path [pattern={Lines[angle=0,distance={4pt}]}] (0.7, -.75) -- (1.1, -.75) -- (1.1, -0.25) -- (0.7, -0.25);
\draw[thin] (0.7, -.75) -- (1.1, -.75) -- (1.1, -0.25) -- (0.7, -0.25) -- (0.7, -.75);
\node[right] at (1.1, 3.5) {$Y \cap \mathbb{Z}^2$};
\node[right] at (1.1, 3) {$[-4, -3) \times [2, \infty)$ for $-\bar{y}^i_1 \leq -4, \bar{y}^i_2 \leq 2$};
\node[right] at (1.1, 2.5) {$[-3, -1) \times [1, 2)$ for $-\bar{y}^i_1 \leq -3, \bar{y}^i_2 \leq 1$};
\node[right] at (1.1, 2) {$[-3, -1) \times [2, 3)$ for $-\bar{y}^i_1 \leq -3, \bar{y}^i_2 \leq 2$};
\node[right] at (1.1, 1.5) {$[-3, -2) \times [3, \infty)$ for $-\bar{y}^i_1 \leq -3, \bar{y}^i_2 \leq 3$};
\node[right] at (1.1, 1) {$[-2, -1) \times [3, \infty)$ for $-\bar{y}^i_1 \leq -2, \bar{y}^i_2 \leq 3$};
\node[right] at (1.1, 0.5) {$[-1, \infty) \times [1, 2)$ for $-\bar{y}^i_1 \leq -1, \bar{y}^i_2 \leq 1$};
\node[right] at (1.1, 0) {$[-1, \infty) \times [2, 3)$ for $-\bar{y}^i_1 \leq -1, \bar{y}^i_2 \leq 2$};
\node[right] at (1.1, -.5) {$[-1, \infty) \times [3, \infty)$ for $-\bar{y}^i_1 \leq -1, \bar{y}^i_2 \leq 3$};
\end{tikzpicture}
\caption{Graph of the sets $V^k$ of Example~\ref{Ex_2}~b)}
\label{Fig_Yk_b}
\vspace{-0.15 in}
\end{figure}
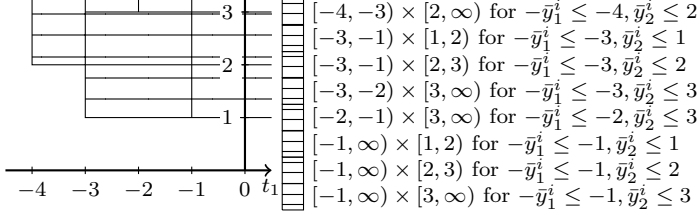

\medskip
\noindent\textit{Example~\ref{Ex_2} (continued) }
For $A = Y^k \subseteq Y \cap \mathbb{Z}^2$ with $A \neq \emptyset$, we get the sets $V^k$ as its graphs are shown in a) Figure~\ref{Fig_Yk_a} and b) Figure~\ref{Fig_Yk_b}.

To conclude this chapter, let us note that based on the sets $V^k$, the functions $\psi$ and $\varphi$ are piecewise constant:
\begin{lemma} \label{Lemma_piecewiseconstant}
Assume \ref{A1}, then the optimal value functions $\psi$ and $\varphi$ are piecewise constant on $\mathrm{dom}\; \Phi$.
\end{lemma}

{\it Proof} For all $k = 1, \ldots, M$, we define \mbox{$\kappa_\psi^k
    := \min_y \left\{d^\top y' \; | \; y' \in Y^k\right\} \in \mathbb{R}$} as well as $\kappa_\varphi^k
    := \min_y \left\{q^\top y \; | \; y \in \mathrm{Argmin}_{y'} \left\{d^\top y' \; | \; y' \in Y^k\right\}\right\} \in \mathbb{R}$. By Lemma~\ref{LemmaA1DomainsPolyhedral} and Lemma~\ref{Lemma_partition}, we have
\[
    \psi(t)
    = \begin{cases} 
        \sum_{k = 1}^M \kappa_\psi^k 1\!\!1_{V^k}(t), &\text{if}\ t \in \mathrm{dom} \; \psi\\
        \infty, &\text{else},
    \end{cases}
\]
as well as
\begin{equation} \label{sumphi}
    \varphi(t)
    = \begin{cases} 
        \sum_{k = 1}^M \kappa_\varphi^k 1\!\!1_{V^k}(t), &\text{if}\ t \in \mathrm{dom} \; \varphi\\
        \infty, &\text{else},
    \end{cases}
\end{equation}
where
$$
    1\!\!1_{V^k}(t) 
    := \begin{cases}
        1, &\text{if}\ t \in V^k \\ 0, &\text{else}
    \end{cases}
$$
denotes the indicator function of the subset $V^k$. \qed

\begin{remark}
Fix $k$ with $V^k \neq \emptyset$. It is $\psi(t') = \psi(t) = \kappa_\psi^k$ and $\varphi(t') = \varphi(t) = \kappa_\varphi^k$ for all $t, t' \in V^k$. Due to Proposition~\ref{Prop_phi}, we have $\varphi(t') = \kappa_\varphi^{k'} \geq \kappa_\varphi^k = \varphi(t)$ for all $t' \in V^{k'}$ and all $t \in V^k$ with $t' \leq t$ and $\psi(t') = \kappa_\psi^{k'} = \kappa_\psi^k = \psi(t)$.
\end{remark}

Alternatively, in the proofs of \cite[Prop. 3.1 3.]{DeMe17} and \cite[Prop. 5.2.3. 3.]{Me17}, the lemma above is proved using the so-called \textit{regions of stability of} $y \in Y \cap \mathbb{Z}^m$: 

$R(y)
    := \left\{t \in \mathrm{dom}\; \varphi \; | \; y \in \Psi(t)\right\}$. 

\noindent The following example shows that the sets $R(y)$ with $y \in Y \cap \mathbb{Z}^m$ are not identical to the sets $V^k$ with $k = 1, \ldots, M$:

\medskip
\noindent\textit{Example~\ref{Ex_2}~(continued) }
Let $d = (1, -1)^\top$. The graph of the sets $R(y)$ are shown in a) Figure~\ref{Fig_Ry_a} and b) Figure~\ref{Fig_Ry_b}.
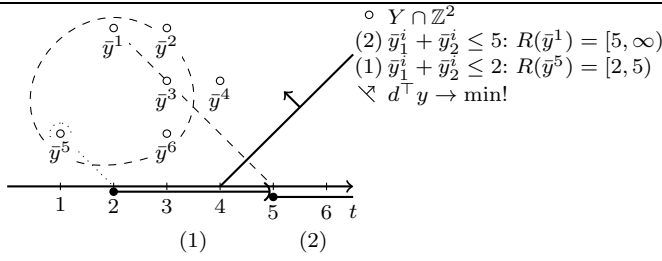
\begin{figure}[ht]
\centering
\vspace{-0.15 in}
\begin{tikzpicture}[scale=0.7]
\draw[thick, ->] (0, 0) -- (6.5, 0);
\node[below] at (6.5, -0.2) {$t$};
\foreach \x in {1,...,6}
    \draw (\x,2pt) -- (\x,-2pt);
\node[below] at (1, -0.05) {$1$};
\node[below] at (2, -0.1) {$2$};
\node[below] at (3, -0.1) {$3$};
\node[below] at (4, -0.1) {$4$};
\node[below] at (5, -0.2) {$5$};
\node[below] at (6, -0.2) {$6$};
\draw[dotted] (1, 1.2) to [closed, curve through =
    {(0.8, 1) (1, 0.8)}] (1.2, 1);
\draw[dashed] (0.6, 1) to [closed, curve through =
    {(1, 0.6) (2.5, 0.5) (3, 0.7) (3.3, 1) (3, 2.8) (2, 3.2) (1, 2.8)}] (0.6, 1);
\draw[dotted] (1, 1) -- (2, 0);
    \filldraw (2, -0.1)circle[radius=2pt];
    \draw[{-)}, thick] (2, -0.1) -- (4.95, -0.1);
\draw[dashed] (2, 3) -- (5, 0);
    \filldraw (5, -0.2)circle[radius=2pt];
    \draw[{-}, thick] (5, -0.2) -- (6.5, -0.2);
\node [rectangle,minimum width=0.4cm,minimum height=0.4cm,fill=white] at (2, 2.65) {};
\node [rectangle,minimum width=0.4cm,minimum height=0.4cm,fill=white] at (3, 2.65) {};
\node [rectangle,minimum width=0.4cm,minimum height=0.4cm,fill=white] at (3, 1.65) {};
\node [rectangle,minimum width=0.4cm,minimum height=0.4cm,fill=white] at (4, 1.65) {};
\node [rectangle,minimum width=0.4cm,minimum height=0.4cm,fill=white] at (1, 0.65) {};
\node [rectangle,minimum width=0.4cm,minimum height=0.4cm,fill=white] at (3, 0.65) {};
\filldraw[white] (1, 1)circle[radius=2pt]
    (3, 1)circle[radius=2pt]
    (3, 2)circle[radius=2pt]
    (4, 2)circle[radius=2pt]
    (2, 3)circle[radius=2pt]
    (3, 3)circle[radius=2pt];
\draw[thin] (1, 1)circle[radius=2pt] node[below] {$\bar{y}^5$}
    (3, 1)circle[radius=2pt] node[below] {$\bar{y}^6$}
    (3, 2)circle[radius=2pt] node[below] {$\bar{y}^3$}
    (4, 2)circle[radius=2pt] node[below] {$\bar{y}^4$}
    (2, 3)circle[radius=2pt] node[below] {$\bar{y}^1$}
    (3, 3)circle[radius=2pt] node[below] {$\bar{y}^2$};
\draw[thick] (4, 0) -- (6.5, 2.5);
\draw[thick, ->] (5.5, 1.5) -- (5.2, 1.8);
\node at (5.75, -1.05) {$(2)$};
\node at (3.5, -1.05) {$(1)$};
\draw[thin] (6.8, 3.25)circle[radius=2pt];
\node at (6.8, 2.75) {$(2)$};
\node at (6.8, 2.25) {$(1)$};
\draw[thin] (6.6, 1.9) -- (6.9, 1.6);
\draw[thin,->] (6.75, 1.75) -- (6.9, 1.9);
\node[right] at (7, 3.25) {$Y \cap \mathbb{Z}^2$};
\node[right] at (7, 2.75) {$\bar{y}^i_1 + \bar{y}^i_2 \leq 5$: $R(\bar{y}^1) = [5, \infty)$};
\node[right] at (7, 2.25) {$\bar{y}^i_1 + \bar{y}^i_2 \leq 2$: $R(\bar{y}^5) = [2, 5)$};
\node[right] at (7, 1.75) {$d^\top y \to \min!$};
\end{tikzpicture}
\caption{Graph of the sets $R(y)$ of Example~\ref{Ex_2}~a)}
\label{Fig_Ry_a}
\vspace{-0.15 in}
\end{figure}
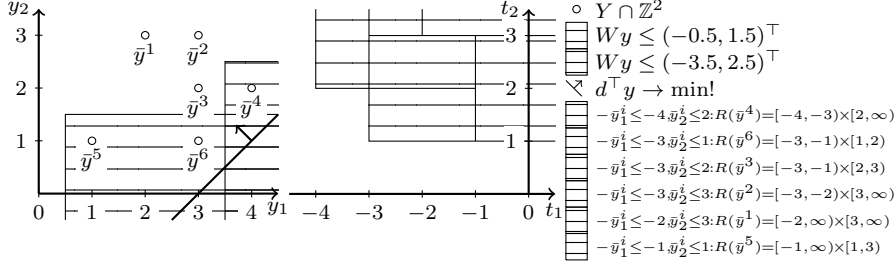
\begin{figure}[ht]
\centering
\begin{tikzpicture}[scale=0.7]
\path [pattern={Lines[angle=67.5,distance={8pt}]}] (4.5, 1.5) -- (0.5, 1.5) -- (0.5, -0.5) -- (4.5, -0.5);
\draw (4.5, 1.5) -- (0.5, 1.5) -- (0.5, -0.5);
\path [pattern={Lines[angle=-67.5,distance={8pt}]}] (4.5, 2.5) -- (3.5, 2.5) -- (3.5, -0.5) -- (4.5, -0.5);
\draw (4.5, 2.5) -- (3.5, 2.5) -- (3.5, -0.5);
\node [rectangle,minimum width=0.4cm,minimum height=0.4cm,fill=white] at (3, 1.65) {};
\node [rectangle,minimum width=0.4cm,minimum height=0.4cm,fill=white] at (4, 1.65) {};
\node [rectangle,minimum width=0.4cm,minimum height=0.4cm,fill=white] at (1, 0.65) {};
\node [rectangle,minimum width=0.4cm,minimum height=0.4cm,fill=white] at (3, 0.65) {};
\node [rectangle,minimum width=0.3cm,minimum height=0.25cm,fill=white] at (4.45, -.3) {};
\draw[thick, ->] (0, 0) -- (4.5, 0) node[below] {$y_1$};
\draw[thick, ->] (0, 0) -- (0, 3.5) node[left] at (0, 3.5) {$y_2$};
\foreach \x in {0,...,4}
    \draw (\x,2pt) -- (\x,-2pt) node[below, fill=white] {$\x$};
\foreach \y in {1,...,3}
    \draw (2pt,\y) -- (-2pt,\y) node[left] {$\y$};
\filldraw[white] (1, 1)circle[radius=2pt]
    (3, 1)circle[radius=2pt]
    (3, 2)circle[radius=2pt]
    (4, 2)circle[radius=2pt]
    (2, 3)circle[radius=2pt]
    (3, 3)circle[radius=2pt];
\draw[thin] (1, 1)circle[radius=2pt] node[below] {$\bar{y}^5$}
    (3, 1)circle[radius=2pt] node[below] {$\bar{y}^6$}
    (3, 2)circle[radius=2pt] node[below] {$\bar{y}^3$}
    (4, 2)circle[radius=2pt] node[below] {$\bar{y}^4$}
    (2, 3)circle[radius=2pt] node[below] {$\bar{y}^1$}
    (3, 3)circle[radius=2pt] node[below] {$\bar{y}^2$};
\draw[thick] (2.5, -.5) -- (4.5, 1.5);
\draw[thick, ->] (4, 1) -- (3.7, 1.3);
\node[white,below] at (4, -.78) {$1$};
\end{tikzpicture}
\hspace*{-1em}
\begin{tikzpicture}[scale=0.7]
\path [pattern={Lines[angle=22.5,distance={8pt}]}] (-4, 2) -- (-3, 2) -- (-3, 3.5) -- (-4, 3.5);
\path [pattern={Lines[angle=45,distance={8pt}]}] (-3, 1) -- (-1, 1) -- (-1, 2) -- (-3, 2);
\path [pattern={Lines[angle=90,distance={8pt}]}] (-3, 2) -- (-1, 2) -- (-1, 3) -- (-3, 3);
\path [pattern={Lines[angle=-45,distance={8pt}]}] (-3, 3) -- (-2, 3) -- (-2, 3.5) -- (-3, 3.5);
\path [pattern={Lines[angle=-22.5,distance={8pt}]}] (-2, 3) -- (0.5, 3) -- (0.5, 3.5) -- (-2, 3.5);
\path [pattern={Lines[angle=0,distance={8pt}]}] (-1, 1) -- (0.5, 1) -- (0.5, 3) -- (-1, 3);
\draw (-4, 3.5) -- (-4, 2) -- (-3, 2) -- (-3, 3.5);
\draw (-3, 2) -- (-3, 1) -- (0.5, 1);
\draw (-3, 2) -- (-1, 2);
\draw (-3, 3) -- (0.5, 3);
\draw (-2, 3) -- (-2, 3.5);
\draw (-1, 1) -- (-1, 3);
\node [rectangle,minimum width=0.25cm,minimum height=0.4cm,fill=white] at (-0.28, 1) {};
\node [rectangle,minimum width=0.28cm,minimum height=0.4cm,fill=white] at (-0.26, 2) {};
\node [rectangle,minimum width=0.28cm,minimum height=0.4cm,fill=white] at (-0.26, 3) {};
\node [rectangle,minimum width=0.28cm,minimum height=0.4cm,fill=white] at (-0.3, 3.5) {};
\draw[thick, ->] (-4.5, 0) -- (0.5, 0) node[below] {$t_1$};
\draw[thick, ->] (0, 0) -- (0, 3.5) node[left] {$t_2$};
\foreach \x in {-4,...,0}
    \draw (\x,2pt) -- (\x,-2pt) node[below] {$\x$};
\foreach \y in {1,...,3}
    \draw (2pt,\y) -- (-2pt,\y) node[left] {$\y$};
\draw[thin] (0.9, 3.5)circle[radius=2pt] node[below] {};
\path [pattern={Lines[angle=67.5,distance={4pt}]}] (0.7, 2.75) -- (1.1, 2.75) -- (1.1, 3.25) -- (0.7, 3.25);
\draw[thin] (0.7, 2.75) -- (1.1, 2.75) -- (1.1, 3.25) -- (0.7, 3.25) -- (0.7, 2.75);
\path [pattern={Lines[angle=-67.5,distance={4pt}]}] (0.7, 2.25) -- (1.1, 2.25) -- (1.1, 2.75) -- (0.7, 2.75);
\draw[thin] (0.7, 2.25) -- (1.1, 2.25) -- (1.1, 2.75) -- (0.7, 2.75) -- (0.7, 2.25);
\draw[thin] (0.7, 2.15) -- (1, 1.85);
\draw[thin,->] (0.85, 2) -- (1, 2.15);
\path [pattern={Lines[angle=22.5,distance={4pt}]}] (0.7, 1.25) -- (1.1, 1.25) -- (1.1, 1.75) -- (0.7, 1.75);
\draw[thin] (0.7, 1.25) -- (1.1, 1.25) -- (1.1, 1.75) -- (0.7, 1.75) -- (0.7, 1.25);
\path [pattern={Lines[angle=45,distance={4pt}]}] (0.7, 0.75) -- (1.1, 0.75) -- (1.1, 1.25) -- (0.7, 1.25);
\draw[thin] (0.7, 0.75) -- (1.1, 0.75) -- (1.1, 1.25) -- (0.7, 1.25) -- (0.7, 0.75);
\path [pattern={Lines[angle=90,distance={4pt}]}] (0.7, 0.25) -- (1.1, 0.25) -- (1.1, 0.75) -- (0.7, 0.75);
\draw[thin] (0.7, 0.25) -- (1.1, 0.25) -- (1.1, 0.75) -- (0.7, 0.75) -- (0.7, 0.25);
\path [pattern={Lines[angle=-45,distance={4pt}]}] (0.7, -.25) -- (1.1, -.25) -- (1.1, 0.25) -- (0.7, 0.25);
\draw[thin] (0.7, -.25) -- (1.1, -.25) -- (1.1, 0.25) -- (0.7, 0.25) -- (0.7, -.25);
\path [pattern={Lines[angle=-22.5,distance={4pt}]}] (0.7, -.75) -- (1.1, -.75) -- (1.1, -.25) -- (0.7, -.25);
\draw[thin] (0.7, -.75) -- (1.1, -.75) -- (1.1, -.25) -- (0.7, -.25) -- (0.7, -0.75);
\path [pattern={Lines[angle=0,distance={4pt}]}] (0.7, -1.25) -- (1.1, -1.25) -- (1.1, -.75) -- (0.7, -.75);
\draw[thin] (0.7, -1.25) -- (1.1, -1.25) -- (1.1, -.75) -- (0.7, -.75) -- (0.7, -1.25);
\node[right] at (1.1, 3.5) {$Y \cap \mathbb{Z}^2$};
\node[right] at (1.1, 3) {$Wy \leq (-0.5, 1.5)^\top$};
\node[right] at (1.1, 2.5) {$Wy \leq (-3.5, 2.5)^\top$};
\node[right] at (1.1, 2) {$d^\top y \to \min!$};
\node[right] at (1.1, 1.5) {\tiny{$-\bar{y}^i_1 \!\!\leq\!\! -4,\!\! \bar{y}^i_2 \!\!\leq\!\! 2\!\!:\!\! R(\bar{y}^4) \!\!=\!\! [-4,\! -3) \!\!\times\!\! [2,\! \infty)$}};
\node[right] at (1.1, 1) {\tiny{$-\bar{y}^i_1 \!\!\leq\!\! -3,\!\! \bar{y}^i_2 \!\!\leq\!\! 1\!\!:\!\! R(\bar{y}^6) \!\!=\!\! [-3,\! -1) \!\!\times\!\! [1,\! 2)$}};
\node[right] at (1.1, 0.5) {\tiny{$-\bar{y}^i_1 \!\!\leq\!\! -3,\!\! \bar{y}^i_2 \!\!\leq\!\! 2\!\!:\!\! R(\bar{y}^3) \!\!=\!\! [-3,\! -1) \!\!\times\!\! [2,\! 3)$}};
\node[right] at (1.1, 0) {\tiny{$-\bar{y}^i_1 \!\!\leq\!\! -3,\!\! \bar{y}^i_2 \!\!\leq\!\! 3\!\!:\!\! R(\bar{y}^2) \!\!=\!\! [-3,\! -2) \!\!\times\!\! [3,\! \infty)$}};
\node[right] at (1.1, -.5) {\tiny{$-\bar{y}^i_1 \!\!\leq\!\! -2,\!\! \bar{y}^i_2 \!\!\leq\!\! 3\!\!:\!\! R(\bar{y}^1) \!\!=\!\! [-2,\! \infty) \!\!\times\!\! [3,\! \infty)$}};
\node[right] at (1.1, -1) {\tiny{$-\bar{y}^i_1 \!\!\leq\!\! -1,\!\! \bar{y}^i_2 \!\!\leq\!\! 1\!\!:\!\! R(\bar{y}^5) \!\!=\!\! [-1,\! \infty) \!\!\times\!\! [1,\! 3)$}};
\end{tikzpicture}
\caption{Graph of a) the minimizers for $t = (-0.5, 1.5)$ and $t = (-3.5, 2.5)$, and b) of the sets $R(y)$ of Example~\ref{Ex_2}~b)}
\label{Fig_Ry_b}
\vspace{-0.15 in}
\end{figure}

\medskip
Now we consider the leader's objective value function $f: \mathbb{R}^n \times \mathbb{R}^s \to \overline{\mathbb{R}}$ based on the notation in \eqref{BSILP}:
\[
    f(x, z)
    := c^\top x + \inf_y \left\{q^\top y \; | \; y \in \Psi(Tx + z)\right\}
    = c^\top x + \varphi(Tx + z)
\]
and continue with a property we need for the proof of Theorem~\ref{Th_stability}:
\begin{corollary}\label{Cor_bounded}
Assume \ref{A1}, then the functions $\psi$ and $\varphi$ are bounded, and there exist $L_f, C_f > 0$ such that $|f(x, z)|
    \leq L_f \|x\| + C_f$ for all $Tx + z \in \mathrm{dom}\; \Phi$.
\end{corollary}

{\it Proof} Based on the considerations in the proof of Lemma~\ref{Lemma_piecewiseconstant} we get 

\noindent $u_\psi 
    := \max_{k = 1, \ldots, M} \kappa^k_\psi
    \geq \psi(t)
    \geq \min_{k = 1, \ldots, M} \kappa^k_\psi
    =: l_\psi$ for all $t \in \mathrm{dom}\; \psi$ as well as $u_\varphi 
    := \max_{k = 1, \ldots, M} \kappa^k_\varphi
    \geq \varphi(t)
    \geq \min_{k = 1, \ldots, M} \kappa^k_\varphi
    =: l_\varphi$ for all $t \in \mathrm{dom}\; \varphi$. Thus, $|\psi(t)|
    \leq \max \left\{|u_\psi|, |l_\psi|\right\}$ and $|\varphi(t)|
    \leq \max \left\{|u_\varphi|, |l_\varphi|\right\}$ for all $t \in \mathrm{dom}\; \Phi$ and $\left|f(x, z)\right|
    \leq \left|c^\top x\right| + \left|\varphi(Tx + z)\right|
    \leq L_f \|x\| + C_f$ for all $Tx + z \in \mathrm{dom}\; \Phi$ with $L_f := \|c\|, \, C_f := \max \left\{|u_\varphi|, |l_\varphi|\right\} \geq 0$.\qed

\section{Two-stage setting} \label{Sec_Model}
Let us now consider a setting in which the parameter $z = Z(\omega)$ in \eqref{BSILP} is the realization of some random vector $Z: \Omega \to \mathbb{R}^s$ on a probability space $(\Omega, \mathcal{F}, \mathbb{P})$, that only the follower can observe. Imposing the information constraint that the leader has to decide in a here-and-now fashion, this gives rise to the bi-level stochastic integer linear program
\begin{equation}\label{BSILP2}
    \min_x \left\{ c^\top x + \inf_y \left\{ q^\top y \; | \; y \in \Psi(Tx + Z(\omega)) \right\} \; | \; x \in X \right\}.
\end{equation}
By definition, the upper level outcome equals $\infty$ whenever the lower level problem has no optimal solution. Deciding nonanticipatorily, the leader thus has to choose $x$ such that $Tx + Z(\omega) \in \mathrm{dom} \; \Psi$ holds for any realization of the randomness. To form a mathematically sound model, we will confine the analysis to situations in which the upper level objective function value is finite as well and thus restrict the leader's choices to the so called induced feasible set given by $F_Z 
    := \{x \; | \; Tx + z \in \mathrm{dom}\;\varphi \  \forall z \in \mathrm{supp} \; \mu_Z\}$ where $\mu_Z := \mathbb{P} \circ Z^{-1} \in \mathcal{P}(\mathbb{R}^s)$ is the Borel probability measure induced by mapping~$Z$.

\begin{example}
In the setting of Example~\ref{ExDomPhi}, set $n = 1$, $T = 1$ and let $\mu_Z$ be given by the Dirac measure at $0$, i.e. $Z(\omega) = 0$ with probability $1$, then $F_Z = \mathrm{dom} \; \varphi = [0,1)$. Hence, $F_Z$ is not closed in general. 
\end{example}

\begin{theorem} \label{Th_F_Z}
Assume \ref{A1} and let $\mathrm{supp} \; \mu_Z$ be polyhedral, then $F_Z$ is polyhedral.
\end{theorem}

{\it Proof} If the support of $\mu_Z$ is given by the union of polyhedra $Q_1, \ldots, Q_k$, we have $F_Z 
    = \bigcap_{i = 1}^k \lbrace x \; | \; Tx + z \in \mathrm{dom}\,\varphi \  \forall z \in Q_i \rbrace$. Since any finite intersection of polyhedral sets is a polyhedral set (cf. last part of this proof), we may restrict the analysis to the case $k = 1$.

Let $\mathrm{supp}\; \mu_Z$ be a nonempty polyhedron $\{z \; | \; Cz \leq c\} =: Q$ and let the domain of function $\varphi$ be a finite union of polyhedra $\mathrm{dom} \; \varphi = \bigcup_{i = 1}^L P_i$ with $P_i := \{t \; | \; A_i t \leq b_i\}$ for matrices $A_i \in \mathbb{R}^{m_i \times s}$ and vectors $b_i \in \mathbb{R}^{m_i}$. For any $x \in F_Z$, we have ($e_{li}(z) := e^\top_l A_i z - e^\top_l b_i + e^\top_l A_i Tx$)
\allowdisplaybreaks
\begin{align}
    & \ \forall\ z \in Q\quad \exists\ i \in \{1, \ldots, L\} \; : \; A_i(Tx + z) \leq b_i \\
    & \Leftrightarrow \forall\ z \in Q\quad \exists\ i \in \{1, \ldots, L\} \quad \forall\ l \in \{1, \ldots, m_i\} \; : \; e^\top_l A_i z \leq e^\top_l b_i - e^\top_l A_i Tx\nonumber\\
    & \Leftrightarrow \sup_{z \in Q} \min_{i \in \{1, \ldots, L\}} \max_{l \in \{1, \ldots, m_i\}} e_{li}(z) \leq 0 \nonumber\\
    & \Leftrightarrow \sup_{z \in Q} \min_\lambda \left\{\sum\nolimits_{i = 1}^L \lambda_i \cdot \max_{l \in \{1, \ldots, m_i\}} e_{li}(z) \; | \; \lambda \in \{0, 1\}^L, \, \sum\nolimits_{i = 1}^L \lambda_i = 1\right\} \leq 0. \nonumber
\end{align}
As the feasible set of the inner minimization problem above is finite, its LP relaxation is solvable and has the same optimal value. Hence,
\allowdisplaybreaks
\begin{align*}
    & \Leftrightarrow\ \sup_{z \in Q} \min_\lambda \left\{\sum\nolimits_{i = 1}^L \lambda_i \cdot \max_{l \in \{1, \ldots, m_i\}} e_{li}(z) \; \Big| \; \lambda \in [0,1]^L, \, \sum\nolimits_{i = 1}^L \lambda_i = 1\right\} \leq 0\\
    & \left.\begin{aligned}
        \Leftrightarrow\ \sup_{z \in Q} \max_{u, w} \bigg\{1\!\!1^\top_L u + w \; \Big| \;
        & u_i + w \leq \max_{l \in \{1, \ldots, m_i\}} e_{li}(z)\ \forall i \in \lbrace 1, \ldots, L\},\\
        & u \leq 0
    \end{aligned}\right\} \leq 0\\
    & \left.\begin{aligned}
        \Leftrightarrow\ \max_{l_1, \ldots, l_L \in \{1, \ldots, s\}} \max_{z, u, w} \Big\{1\!\!1^\top_L u + w \; \big| \;
        & e^\top_{l_i} u + w \leq e_{l_ii}(z)\ \forall \; i \in \{1, \ldots, L \},\\
        & u \leq 0, \, Cz \leq c
    \end{aligned}\right\} \leq 0
\end{align*}    
Since setting $u = 0$, $z = z_0$ for some fixed $z_0 \in \mathrm{supp}\; \mu_Z$, and
$$
    w 
    = \min \left\{1, \min_{i \in \{1, \ldots, L\}} e_{l_ii}(z_0)\right\}
$$
yields a feasible point for the modified problem, we may add the restriction $1\!\!1^\top_L u + w \leq 1$ to ensure that the inner maximization problem has an optimal solution. By linear programming duality, we may continue the previous equivalences:
\allowdisplaybreaks
\begin{align*}
    & \left.\begin{aligned}
        \Leftrightarrow\ 
        \max_{l_1, \ldots, l_L \in \{1, \ldots, s\}} \max_{z, u, w} \big\{ & 1\!\!1^\top_L u + w \; \big|\\
        & 1\!\!1^\top_L u + w \leq 1, \, u \leq 0, \, Cz \leq c,\\
        & e^\top_{l_i} u + w \leq e_{li}(z) \ \forall \; i \in \{1, \ldots, L \} 
    \end{aligned}\right\} \leq 0\\
    & \left.\begin{aligned}
        \Leftrightarrow\ 
        \max_{l_1, \ldots, l_L \in \{1, \ldots, s\}} \min_{\alpha, \beta, \gamma} \bigg\{
        & \sum\nolimits_{i = 1}^L \alpha_i \left(-e^\top_{l_i} b_i + e^\top_{l_i} A_i Tx\right) + \beta^\top c + \gamma \; \big|\\
        & 1\!\!1_L \gamma + \sum\nolimits_{i = 1}^L e_{l_i} \alpha_i \leq 1\!\!1_L, \, \sum\nolimits_{i = 1}^L \alpha_i + \gamma = 1,\\
        & \sum\nolimits_{i = 1}^L (-A_i^\top e_{l_i})\alpha_i + C^\top \beta = 0, \, \alpha, \beta, \gamma \geq 0
    \end{aligned}\right\} \leq 0.
\end{align*}
Denoting the finite nonempty set of vertices of the feasible set of the inner minimization problem above by $\mathcal{V}(l_1, \ldots, l_L)$ and the $i$-th biggest element in $\{1, \ldots, s\}^L$ with respect to the lexicographical order by $\mathbb{L}_i$, we may reformulate the previous condition as
\allowdisplaybreaks
\begin{align*}
    & \ x \in \bigcap_{\substack{l_1, \ldots, l_L\\\in \{1, \dots, s\}}} \bigcup_{\substack{(\alpha, \beta, \gamma)\\\in \mathcal{V}(l_1, \ldots, l_L)}} \underbrace{\left\{x' \; \Big| \; \sum\nolimits_{i = 1}^L \alpha_i \left(-e^\top_{l_i} b_i + e^\top_{l_i} A_i Tx'\right) + \beta^\top c + \gamma \leq 0\right\}}_{=: P(l_1, \ldots, l_L, \alpha, \beta, \gamma)}\\
    & \Leftrightarrow x \in \bigcup_{(\rho_1, \ldots, \rho_{s^L}) \in \prod_{i = 1}^{s^L} \mathcal{V}(\mathbb{L}_i)} \bigcap_{j = 1}^{s^L} P(\mathbb{L}_j, \rho_j),
\end{align*}
which yields the desired representation. \qed

\begin{lemma} \label{Lemma_welldefined}
Assume \ref{A1}, then the function $\mathbb{F}: F_Z \to \mathcal{L}^\infty(\Omega, \mathcal{F}, \mathbb{P})$ with $\mathbb{F}(x)(\cdot) := \varphi(Tx + Z(\cdot))$ is well-defined.
\end{lemma}

{\it Proof} As $\varphi$ is measurable (cf. Lemma~\ref{Lemma_phimeasurable2}), we conclude that $\mathbb{F}(x)(\cdot)$ is measurable, i.e. $\mathbb{F} \in \mathcal{L}^0(\Omega, \mathcal{F}, \mathbb{P})$. Moreover, \ref{A1} implies
$$ 
    \|\mathbb{F}(x) \|_{\mathcal{L}^\infty(\Omega, \mathcal{F}, \mathbb{P})} 
    \leq \max_y \lbrace |q^\top y| \; | \; y \in Y \cap \mathbb{Z}^m\rbrace 
    < +\infty
$$
and thus $\mathbb{F} \in \mathcal{L}^\infty(\Omega, \mathcal{F}, \mathbb{P})$. \qed

\begin{theorem} \label{Th_FContinuous2}
Assume \ref{A1}, $\mu_Z\left[D_\varphi(x)\right] = 0$, where $D_\varphi: F_Z \rightrightarrows \mathbb{R}^s$ with $D_\varphi(x) 
    := \left\{z \; | \; \varphi\ \text{is discontinuous at}\ (Tx + z)\right\}$, then the mapping $\mathbb{F}$ is continuous at $x \in F_Z$ with respect to any $\mathcal{L}^p$-norm with $p \in [1, \infty)$.
\end{theorem}

{\it Proof} Let $x \in F_Z$ and $\{x_l\}_{l \in \mathbb{N}} \subseteq F_Z$ any sequence with $x_l \to x$ for $l \to \infty$. With the reformulation of $\varphi$ at \eqref{sumphi} and
\begin{align*}
    |\varphi(Tx_l + z)|
    = \left|\sum\nolimits_{k = 1}^M \kappa_\varphi^k 1\!\!1_{V^k}(Tx_l + z) \right|
    \leq \underbrace{\max_{k = 1, \ldots, M} |\kappa_\varphi^k|}_{=: \bar{\kappa}}
    < \infty
\end{align*}
$\mu_Z$-almost everywhere for all $l \in \mathbb{N}$ we have a measurable majorant $\bar{\kappa}$. Based on $\mu_Z\left[D_\varphi(x)\right] = 0$, we have $\lim_{l \to \infty} \varphi(Tx_l + z) 
    = \varphi(Tx + z)$ $\mu_Z\text{-almost everywhere}$ for the functions $\mathbb{F}(x_l) \in \mathcal{L}^p(\Omega, \mathcal{F}, \mathbb{P})$ with $l \in \mathbb{N}$. We obtain \mbox{$\mathbb{F}(x) \in \mathcal{L}^p(\Omega, \mathcal{F}, \mathbb{P})$} as well as
\begin{align*}
    & \lim_{l \to \infty} \left\|\mathbb{F}(x_l) - \mathbb{F}(x)\right\|_{\mathcal{L}^p(\Omega, \mathcal{F}, \mathbb{P})}\\
    & = \lim_{l \to \infty} \left(\int_{\mathbb{R}^s} \left|\varphi(Tx_l + z) - \varphi(Tx + z)\right|^p ~\mu_Z(\mathrm{d}z)\right)^{1/p}
    = 0
\end{align*}
using a majorized convergence theorem for $\mathcal{L}^p$-functions, cf. \cite[\S 12 Th. 4]{Fo12}. \qed

\smallskip
A characterisation of the sets $D_\varphi$ is given in the proof of the next corollary.
\begin{corollary} \label{Cor_FContinuous}
Assume \ref{A1} and that the Borel measure $\mu_Z$ is absolutely continuous with respect to the Lebesgue measure, then the mapping $\mathbb{F}$ is continuous with respect to any $\mathcal{L}^p$-norm with $p \in [1, \infty)$.
\end{corollary}

{\it Proof} By \cite[Prop. 1.1 (i)]{GoJoRo03}, we have
\[
    \bar{V}^k
    := \mathrm{cl} \; V^k
    = \bigcap_{y \in Y^k} \{t \; | \; Wy \leq t\} \cap
    \bigcap_{\hat{y} \in \left(Y \cap \mathbb{Z}^m\right) \setminus Y^k} \bigcup_{j = 1}^s \left\{t \; | \; e_j^\top W\hat{y} \geq t_j\right\}
\]
and thus $\bar{V}^k
    \subseteq V^k \cup \bigcup_{\hat{y} \in \left(Y \cap \mathbb{Z}^m\right) \setminus Y^k} \bigcup_{j = 1}^s \left\{t \; | \; e_j^\top W\hat{y} = t_j\right\}$. In particular, the set \mbox{$\bar{V}^k \setminus V^k$} is contained in a finite union of hyperplanes and therefore a Lebesgue null set. We also obtain the still missing characterization of the set of discontinuity points from Theorem~\ref{Th_FContinuous2}: 
\begin{equation}
    \left(\bar{V}^k \setminus V^k\right) \oplus (-Tx) \supseteq D_\varphi(x)
    \quad \text{for all}\ x \in F_Z.\label{nullset}
\end{equation}
Based on that, the image of the setvalued mapping $D_\varphi$ is a Lebesgue null set. Due to the absolute continuity of $\mu_Z$ with respect to the Lebesgue measure, we get continuity of the function $\mathbb{F}$ at all $x \in F_Z$ with respect to any $\mathcal{L}^p$-norm with $p \in [1, \infty)$ by the help of Theorem~\ref{Th_FContinuous2}. \qed

\smallskip
The following example shows that the assumptions of Corollary~\ref{Cor_FContinuous} are too weak to guarantee the local Lipschitz continuity of $\mathbb{F}$, even if the support of $\mu_Z$ is compact:
\begin{example} \label{ExLocLipschitz}
Consider the case where $n = m = s = 1$, $W = T = q = 1$, $d = -1$,
and $Y = [0,1]$, cf. Figure~\ref{Fig_LocLipschitz}, then
\begin{figure}[ht]
\centering
\vspace{-0.15 in}
\begin{tikzpicture}[scale=0.7]
\path [pattern={Lines[angle=22.5,distance={8pt}]}] (-1.5, -0.5) -- (0.7, -0.5) -- (0.7, 0.3) -- (-1.5, 0.3);
\draw (0.7, -0.75) -- (0.7, 0.3);
\foreach \x in {-1,...,2}
  \draw (\x,2pt) -- (\x,-2pt) node[below,fill=white] {$\x$};
\draw[thick, ->] (-1.5, 0) -- (2.5, 0) node[below] {$y$};
\draw[thick, ->] (-1.5, -0.6) -- (2.5, -0.6) node[below] {$t$};
\filldraw (0, -0.65)circle[radius=2pt];
    \draw[thick] (0, -0.65) -- (0.7, -0.65);
    \node[below] at (0.35, -0.65) {$t \in \mathrm{dom}\; \varphi$};
\filldraw[white] (0, 0)circle[radius=2pt]
(1, 0)circle[radius=2pt];
\draw[thin] (0, 0)circle[radius=2pt]
(1, 0)circle[radius=2pt];
\draw[thick] (-0.5, -0.5) -- (-0.5, 0.3);
\draw[thick, ->] (-0.5, -0.1) -- (-0.2, -0.1);
\end{tikzpicture}
\begin{tikzpicture}[scale=0.7]
\path [pattern={Lines[angle=22.5,distance={8pt}]}] (-1.5, -0.5) -- (1.7, -0.5) -- (1.7, 0.5) -- (-1.5, 0.5);
\draw (1.7, -0.75) -- (1.7, 0.5);
\foreach \x in {-1,...,2}
  \draw (\x,2pt) -- (\x,-2pt) node[below,fill=white] {$\x$};
\draw[thick, ->] (-1.5, 0) -- (2.5, 0) node[below] {$y$};
\draw[thick, ->] (-1.5, -0.6) -- (2.5, -0.6) node[below] {$t$};
\filldraw (0, -0.65)circle[radius=2pt];
    \draw[{-)},thick] (0, -0.65) -- (1, -0.65);
    \node[below] at (0.5, -0.65) {$t \in \mathrm{dom}\; \varphi$};
\filldraw[white] (0, 0)circle[radius=2pt]
(1, 0)circle[radius=2pt];
\draw[thin] (0, 0)circle[radius=2pt]
(1, 0)circle[radius=2pt];
\draw[thick] (0.5, -0.5) -- (0.5, 0.3);
\draw[thick, ->] (0.5, -0.1) -- (0.8, -0.1);
\draw[thin] (2.9, 0.25)circle[radius=2pt];
\path [pattern={Lines[angle=22.5,distance={4pt}]}] (2.7, -.5) -- (3.1, -.5) -- (3.1, 0) -- (2.7, 0);
\draw[thin] (2.7, -.5) -- (3.1, -.5) -- (3.1, 0) -- (2.7, 0) -- (2.7, -.5);
\draw[thin] (2.7, -.6) -- (3, -.9);
\draw[thin,->] (2.85, -.75) -- (3, -.6);
\node[right] at (3.1, 0.25) {$Y \cap \mathbb{Z}^2$};
\node[right] at (3.1, -0.25) {$Wy \leq t$};
\node[right] at (3.1, -0.75) {$d^\top y \to \min!$};
\end{tikzpicture}
\caption{Graph of the feasible points and $t \in \mathrm{dom}\; \varphi$ with a) $t \in (0, 1)$ and b) $t > 1$ of Example~\ref{ExLocLipschitz}}
\label{Fig_LocLipschitz}
\vspace{-0.15 in}
\end{figure}
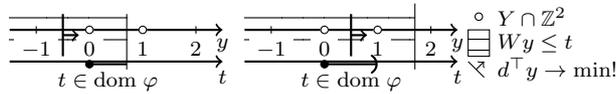
$$
    \varphi(t) 
    = \begin{cases} \infty, &\text{if} \; t < 0 \\
    0, &\text{if} \; t \in [0,1) \\
    1, &\text{if} \; t \geq 1. \end{cases}
$$
Moreover, let $\mu_Z \in \mathcal{P}(\mathbb{R})$ be the Borel probability measure with Lebesgue density
$$
    \delta_Z(z) 
    := \begin{cases} \frac{1}{2\sqrt{z}}, &\text{if} \; z \in (0,1] \\
    0, &\text{else}, \end{cases}
$$
then
$$    \|\mathbb{F}(x)\|_{\mathcal{L}^1(\Omega, \mathcal{F}, \mathbb{P})} 
    = \int\limits_{0}^1 \frac{1}{2\sqrt{z}} \varphi(x+z)~\mathrm{d}z 
    = \int\limits_{1-x}^1 \frac{1}{2\sqrt{z}}~\mathrm{d}z 
    = 1 - \sqrt{1-x}
$$
holds for any $x \in [0,1]$. Hence, $\mathbb{F}$ is not locally Lipschitz continuous at $1$ w.r.t. any $\mathcal{L}^p$-norm with $p \in [1, \infty]$.
\end{example}

Furthermore, Corollary~\ref{Cor_FContinuous} cannot be extended to $p = \infty$:
\begin{example} 
Consider the setting of Example~\ref{ExLocLipschitz}, but let $Z$ be uniformly distributed on $[0,1]$, then
$$
\mathbb{F}(x)(\omega) := \begin{cases} 0, &\text{if} \; x + Z(\omega) < 1 \\ 1, &\text{else} \end{cases}
$$
holds for any $x \in F_Z = [0,\infty)$. Thus, $\mathbb{F}$ is not continuous at $0$ w.r.t. the  $\mathcal{L}^\infty$-norm, as it holds
\begin{align*}
    &\left\|\mathbb{F}\left(\tfrac{1}{j}\right) - \mathbb{F}(0)\right\|_{\mathcal{L}^\infty(\Omega, \mathcal{F}, \mathbb{P})} 
    = \underset{\omega \in \Omega}{\mathrm{ess} \sup} \left|\mathbb{F}\left(\tfrac{1}{j}\right)(\omega) - \mathbb{F}(0)(\omega)\right|\\
    &= \underset{\omega \in \Omega}{\mathrm{ess} \sup} \left|\varphi\left(\tfrac{1}{j} + Z(\omega)\right) - \varphi(Z(\omega))\right|
    = 1
    \quad\text{for any}\ j \in \mathbb{N}.
\end{align*}
\end{example}

The subsequent analysis relies on the following additional assumption:
\begin{enumerate}[label=$\mathrm{(A\arabic*)}$,leftmargin=1cm,start=2]
    \item $\mu_Z$ has a uniformly bounded Lebesgue density.\label{A2}
\end{enumerate}
\begin{theorem} \label{ThFLipschitz}
Assume \ref{A1}, \ref{A2} as well as $p \in [1, \infty)$ and let the support of $\mu_Z$ be compact, then the mapping $\mathbb{F}$ is Hölder continuous with exponent $\frac{1}{p}$ with respect to the $\mathcal{L}^p$-norm.
\end{theorem}

{\it Proof} Based on the notation in the proof of Corollary~\ref{Cor_FContinuous}, we introduce the set-valued mapping $\Delta^{k,j}: F_Z \times F_Z \rightrightarrows \mathbb{R}^s$ defined by 

$\Delta^{k,j}(x, x') 
    := \big(\bar{V}^k \oplus (- Tx)\big) \cap \big(\bar{V}^j \oplus (- Tx')\big)$,
    
\noindent and fix some $k, j \in \{1, \ldots, M\}$ and some $x \in F_Z$.
As $\bar{V}^k$ and $\bar{V}^j$ are closed, $\Delta^{k,j}(x, \cdot)$ is outer semicontinuous (cf. \cite[Lemma A.3.]{BuClDe20}), i.e., 

$\lim \sup_{x' \to x_0} \Delta^{k, j}(x, x') \subseteq \Delta^{k, j}(x, x_0)$ for all $x_0 \in F_Z$.

For $j \neq j'$ the sets $V^j$ and $V^{j'}$ are disjoint, cf. Lemma~\ref{Lemma_partition}, and the set $\big(\bar{V}^j \oplus (- Tx')\big) \cap \big(\bar{V}^{j'} \oplus (- Tx')\big)$ is a null set. Due to Lemma~\ref{LemmaA1DomainsPolyhedral} and Lemma~\ref{Lemma_partition}, $\bigcup_{j = 1}^M \bar{V}^j
    = \bigcup_{j = 1}^M V^k
    = \mathrm{dom}\, \varphi$ and it follows
    
    $\sum_{j = 1}^M \mu_Z\left[\Delta^{k,j}(x, x')\right]
    = \mu_Z\left[\bar{V}^k \oplus (- Tx)\right]$.
    
\noindent Invoking \cite[Lemma 3.4]{Cl21}, we obtain continuity of a mapping \mbox{$M^{k,j}: F_Z^2 \to [0, 1]$} on $\{x\} \times F_Z$ defined by \mbox{$M^{k,j}(x, x') := \mu_Z\left[\Delta^{k,j}(x, x')\right]$} for any $j \in \{1, \ldots, M\}$.

Fix an arbitrary point $x \in F_Z$ and consider any $x' \in F_Z$. We can write
\allowdisplaybreaks
\begin{align}
    & \|\mathbb{F}(x) - \mathbb{F}(x')\|^p_{\mathcal{L}^p(\Omega, \mathcal{F}, \mathbb{P})}\\
    & = \int_{\mathbb{R}^s} \left|\sum\nolimits_{k = 1}^M \kappa_\varphi^k 1\!\!1_{V^k}(Tx + z) - \sum\nolimits_{j = 1}^M \kappa_\varphi^j 1\!\!1_{V^j}(Tx' + z)\right|^p ~\mu_Z(\mathrm{d}z) \nonumber\\
    & = \sum\nolimits_{k = 1}^M \sum\nolimits_{j = 1}^M \int_{\big(V^k \oplus (- Tx)\big) \cap \big(V^j \oplus (- Tx')\big)} |\kappa_\varphi^k - \kappa_\varphi^j|^p ~\mu_Z(\mathrm{d}z) \nonumber\\
    & \leq \sum\nolimits_{k = 1}^M \sum\nolimits_{\substack{j = 1\\k \neq j}}^M \underbrace{\max_{k, j = 1, \ldots, M} |\kappa_\varphi^k - \kappa_\varphi^j|^p}_{=: \alpha < +\infty} \mu_Z\left[\big(V^k \oplus (- Tx)\big) \cap \big(V^j \oplus (- Tx')\big)\right] \nonumber\\
    & \leq \sum\nolimits_{k = 1}^M \sum\nolimits_{\substack{j = 1\\k \neq j}}^M \alpha M^{k, j}(x, x') \label{sumMkj}
\end{align}
based on the absolute continuity of $\mu_Z$. Now it suffices to derive Lipschitz estimates for $M^{k, j}(x, x')$ for all $k, j = 1, \ldots, M$ and for all $x \in F_Z$ to continue estimation \eqref{sumMkj}:
\begin{align}
    & M^{k, j}(x, x')
    = \mu_Z\left[\big(\bar{V}^k \oplus (- Tx)\big) \cap \big(\bar{V}^j \oplus (- Tx')\big)\right] \nonumber\\
    & = \mu_Z\left[\big(\bar{V}^k \oplus (- Tx)\big) \cap \big(\bar{V}^j \oplus (- Tx + T(x - x'))\big)\right] \label{Mijx}
\end{align}
For $k \neq j$, it is $V^k \cap V^j = \emptyset$ and we get
\begin{itemize}
    \item $\big(\bar{V}^k \oplus (- Tx)\big) \cap \big(\bar{V}^j \oplus (- Tx)\big)$ is null set and, in particular, it is contained in a finite union of hyperplanes, and
    \item $\big(\bar{V}^k \oplus (- Tx)\big) \cap \big(\bar{V}^j \oplus (- Tx)\big) 
    \subseteq \mathrm{bd}\, \big(\bar{V}^j \oplus (- Tx)\big)
    = (\mathrm{bd}\,\bar{V}^j) \oplus (- Tx)$
\end{itemize}
such that we receive
\begin{align*}
    & \big(\bar{V}^k \oplus (- Tx)\big) \cap \big(\bar{V}^j \oplus (- Tx + T(x - x'))\big)\\
    & \subseteq \mathrm{bd}\, \bar{V}^j \oplus (- Tx) \oplus [0, 1]\, T(x - x')\\
    & = \left\{v + lT(x - x') \; | \; v \in \mathrm{bd}\, \bar{V}^j \oplus (- Tx), \, l \in [0, 1]\right\}.
\end{align*}
Now we can proceed with the estimation of \eqref{Mijx} by
\[
    \leq \delta \lambda^s\left[\left\{v + l T(x - x') \; | \; v \in \mathrm{bd}\, \bar{V}^j \oplus (- Tx), \, l \in [0, 1]\right\} \cap \mathrm{supp}\, \mu_Z\right],
\]
where $\delta < \infty$ is a bound for the Lebesgue density of $\mu_Z$. The boundary of $\bar{V}^j \oplus (- Tx)$ is the finite union of $K$ (independent of $k, j$, $x, x'$) hyperplanes $H_1, \ldots, H_K$. By the help of Cavalieri's principle and $\mathrm{diam}\,(\mathrm{supp}\,\mu_Z)^{s - 1} = \beta$, we finish the proof
\begin{align*}
    \leq \delta \sum\nolimits_{l = 1}^K \lambda^{s - 1}\left[H_l \cap \mathrm{supp}\, \mu_Z\right] \cdot \|T(x - x')\|
    \leq \delta \cdot K \beta \cdot \|T\| \cdot \|x - x'\|.\tag*{$\qed$}
\end{align*}

Example~\ref{ExLocLipschitz} shows that boundedness of the Lebesgue density is essential for Theorem~\ref{ThFLipschitz}. If the support of $\mu_Z$ is not bounded, we still obtain a weaker estimate:
\begin{theorem}
Assume \ref{A1}, \ref{A2} and $p \in [1, \infty)$, then for any $\epsilon > 0$ there exists a real number $L(\epsilon)$ such that $\|\mathbb{F}(x) - \mathbb{F}(x')\|_{\mathcal{L}^p(\Omega, \mathcal{F}, \mathbb{P})} \leq L(\epsilon) \|x - x\|^{\frac{1}{p}}$ holds for any $x, x' \in F_Z$.
\end{theorem}

{\it Proof} Throughout the proof, we will fix some $\epsilon > 0$ and use the notation established in the proof of Theorem~\ref{ThFLipschitz}. By construction, we may assume that $K \geq 2$. As the probability measure $\mu_Z$ is tight by \cite[Th. 1.3]{Bi99}, there is some compact set $C(\epsilon) \subset \mathbb{R}^s$ such that
$$
\mu_Z \left[ \mathbb{R}^s \setminus C(\epsilon) \right] \leq \frac{\epsilon^p}{K(K-1)\alpha}.
$$
By the arguments given in the proof of Theorem~\ref{ThFLipschitz}, we thus have
\allowdisplaybreaks
\begin{align*}
    & M^{k, j}(x, x')\\
    & \leq \mu_Z\left[\left\{v + l T(x - x') \; | \; v \in \mathrm{bd}\, \bar{V}^j \oplus (- Tx), \, l \in [0, 1]\right\} \cap C(\epsilon) \right] + \tfrac{\epsilon^p}{K(K-1)\alpha}\\
    & \leq \delta K \mathrm{diam}(C(\epsilon))^{s-1} \|T\|\cdot \|x - x'\| + \tfrac{\epsilon^p}{K(K-1)\alpha}
\end{align*}
and \eqref{sumMkj} yields the desired estimate
\allowdisplaybreaks
\begin{align*}
    &\|\mathbb{F}(x) - \mathbb{F}(x')\|_{\mathcal{L}^p(\Omega, \mathcal{F}, \mathbb{P})} 
    \leq \left(\delta K \mathrm{diam}(C(\epsilon))^{s-1} \|T\|\cdot \|x - x'\| + \epsilon^p \right)^{\frac{1}{p}}\\
    & \leq \left( \delta K \mathrm{diam}(C(\epsilon))^{s-1} \|T\| \right)^{\frac{1}{p}}  \|x - x'\|^{\frac{1}{p}} + \epsilon,
\end{align*}
where the second inequality follows from the fact that the function $t \to t^{\frac{1}{p}}$ is subadditive on $[0,\infty)$. \qed

\section{Risk-averse approach} \label{Sec_4}
As a first choice, we might evaluate the random upper level objective function based on its expected value, i.e. consider the risk neutral bi-level stochastic program $\min_{x \in F_Z} \left\{\mathbb{E} \left[\mathbb{F}(x)\right]\right\}$, which is well-defined by Lemma~\ref{Lemma_welldefined}. More in general, to allow for varying degrees of risk aversion, we take into account a mapping $\mathcal{R}\colon \mathcal{X} \to \mathbb{R}$ with $\mathcal{L}^\infty(\Omega, \mathcal{F}, \mathbb{P}) \subseteq \mathcal{X} \subseteq \mathcal{L}^0(\Omega, \mathcal{F}, \mathbb{P})$ and consider the bi-level stochastic program
\begin{equation} \label{StochasticBi-levelProgram}
    \min_{x \in F_Z} \left\{\mathcal{R}\left[\mathbb{F}(x)\right]\right\}.
\end{equation}
$\mathcal{R}$ will typically be a monetary risk measure in the sense of \cite[Def. 4.1]{FoSc11} meaning it satisfies the following conditions:
\begin{itemize}
    \item Monotonicity: $\mathcal{R}[Y_1] \leq \mathcal{R}[Y_2]$ for all $Y_1, Y_2 \in \mathcal{X}$ satisfying $Y_1 \leq Y_2$ $\mathbb{P}$-almost surely.
    \item Translation equivariance: $\mathcal{R}[Y + m] = \mathcal{R}[Y] + m$ for all $Y \in \mathcal{X}$  and $m \in \mathbb{R}$.
\end{itemize}
Moreover, we will assume the following:
\begin{enumerate}[label=$\mathrm{(A\arabic*)}$,leftmargin=1cm,start=3]
    \item $\mathcal{R}\colon \mathcal{L}^p(\Omega, \mathcal{F}, \mathbb{P}) \to \mathbb{R}$ with some $p \in [1, \infty)$ is convex and nondecreasing as defined above.\label{A3}
\end{enumerate}
\begin{remark}
\ref{A3} holds for any convex risk measure in the sense of \cite{FrRoGi02} and \cite{FoSc02}, i.e. for any monetary risk measure that is convex. In particular, this includes the expectation, the mean-upper semideviation of any order and the Conditional Value-at-Risk. However, as we do not assume translation equivariance, the assumption is also fulfilled for the expected excess of arbitrary order (cf. \cite[\S 6]{ShDeRu09}).
\end{remark}

\begin{theorem} \label{Th_QContinuous}
Assume \ref{A1}, \ref{A3}, then the following statements hold true:
\begin{enumerate}
    \item The function $\mathcal{Q}_\mathcal{R} \colon F_Z \to \mathbb{R}$ defined by $\mathcal{Q}_\mathcal{R}(x) := \mathcal{R}\left[\mathbb{F}(x)\right]$ is continuous at $x \in F_Z$ if $\mu_Z\left[D_\varphi(x)\right] = 0$.
    \item The function $\mathcal{Q}_\mathcal{R}$ is continuous if the Borel measure $\mu_Z$ is absolutely continuous with respect to the Lebesgue measure. In particular, the bi-level stochastic problem \eqref{StochasticBi-levelProgram} has an optimal solution whenever the induced feasible set $F_Z$ is nonempty and compact.
\end{enumerate}
\end{theorem}

{\it Proof} As $\mathcal{R}$ is continuous by \cite[Lemma 3]{BuClCoRuSaSc21}, the results follow from Lemma~\ref{Lemma_welldefined} and Theorem~\ref{Th_FContinuous2} or Corollary~\ref{Cor_FContinuous}. \qed

\begin{remark}
In general, assumption \ref{A3} does not hold for a so-called certainty equivalent $\mathrm{CE}_u$ (cf. \cite[Def. 2.8]{BuClDe20}) because of the lack of convexity. However, Lipschitz continuity is guaranteed under the assumption, that $u: \mathbb{R} \to \mathbb{R}$ is strictly increasing and bi-Lipschitz, cf. \cite[Lemma 3.6]{BuClDe20}.
\end{remark}

The following result for bi-level linear optimization under uncertainty can be found in \cite[Prop. 1]{BuCl20}:
\begin{proposition}  \label{PropQRLipschitz}
Assume \ref{A1}, \ref{A2} and let the support of $\mu_Z$ be compact. Then the following statements hold true for any mapping $\mathcal{R}: L^p(\Omega, \mathcal{F}, \mathbb{P}) \to \mathbb{R}$:
\begin{enumerate}
    \item $\mathcal{Q}_\mathcal{R}$ is locally Hölder continuous with exponent $\frac{1}{p}$ if $\mathcal{R}$ is convex and continuous.
    \item $\mathcal{Q}_\mathcal{R}$ is locally Hölder continuous with exponent $\frac{1}{p}$ if we assume \ref{A3}.
    \item $\mathcal{Q}_\mathcal{R}$ is locally Hölder continuous with exponent $\frac{1}{p}$ if $\mathcal{R}$ is a convex risk measure.
    \item $\mathcal{Q}_\mathcal{R}$ is Hölder continuous with exponent $\frac{1}{p}$ if $\mathcal{R}$ is Lipschitz continuous.
    \item $\mathcal{Q}_\mathcal{R}$ is Hölder continuous with exponent $\frac{1}{p}$ if $\mathcal{R}$ is a coherent risk measure.
\end{enumerate}
\end{proposition}

{\it Proof} See the proof of \cite[Prop. 1]{BuCl20} in combination with Theorem~\ref{ThFLipschitz}. \qed

\begin{remark}
We obtain (local) Lipschitz continuity for all feasible risk measures with $p = 1$, which includes, for example, the expectation, the expected excess of order $1$, and the Conditional Value-at-Risk.
\end{remark}

Due to the lack of convexity, Theorem~\ref{Th_QContinuous} and the subsequent proposition do not apply to the excess probability and the Value-at-Risk. However, the arguments from the proof of Theorem~\ref{Th_QContinuous} can be used another time:
\begin{theorem}
Assume \ref{A1} and that the Borel measure $\mu_Z$ is absolutely continuous with respect to the Lebesgue measure. Fix $\eta \in \mathbb{R}$, then the function $\mathcal{Q}_{\mathrm{EP}_\eta}(x) := \mathrm{EP}_\eta\left[\mathbb{F}(x)\right] = \mu_Z\left[\left\{\omega \in \Omega \; | \; \mathbb{F}(x)(\omega) > \eta\right\}\right]$ is continuous. Furthermore, let the induced feasible set $F_Z$ be nonempty and compact. Then $\min_{x \in F_Z} \left\{\mathcal{Q}_{\mathrm{EP}_\eta}(x)\right\}$ is solvable.
\end{theorem}

{\it Proof} The function $\mathcal{Q}_{\mathrm{EP}_\eta}(x)$ is real-valued on $F_Z$ due to Lemma~\ref{Lemma_welldefined}. Fix an arbitrary point $x \in F_Z$ and consider any $x' \in F_Z$. By Lemma~\ref{Lemma_partition}, we have
\allowdisplaybreaks
\begin{align*}
    &\left|\mathcal{Q}_{\mathrm{EP}_\eta}(x) - \mathcal{Q}_{\mathrm{EP}_\eta}(x')\right|\\
    & = \left|\mu_Z\left[\left\{\omega \in \Omega \; | \; \varphi(Tx + Z(\omega)) > \eta\right\}\right] - \mu_Z\left[\left\{\omega \in \Omega \; | \; \varphi(Tx' + Z(\omega)) > \eta\right\}\right]\right|\\
    & = \left|\mu_Z\left[\left\{\omega \in \Omega \; | \; \sum\nolimits_{k = 1}^M \kappa_\varphi^k 1\!\!1_{V^k}(Tx + Z(\omega)) > \eta\right\}\right]\right.\\
    &\quad \left.- \mu_Z\left[\left\{\omega \in \Omega \; | \; \sum\nolimits_{j = 1}^M \kappa_\varphi^j 1\!\!1_{V^j}(Tx' + Z(\omega)) > \eta\right\}\right]\right|\\
    &= \left|\mu_Z\left[\bigcup\nolimits_{\substack{k = 1\\\kappa_\varphi^k > \eta}}^M V^k \oplus (-Tx)\right] - \mu_Z\left[\bigcup\nolimits_{\substack{j = 1\\\kappa_\varphi^j > \eta}}^M V^j \oplus (-Tx')\right]\right|\\
    & = \left|\sum\nolimits_{\substack{k = 1\\\kappa_\varphi^k > \eta}}^M \mu_Z\left[V^k \oplus (-Tx)\right] - \sum\nolimits_{\substack{j = 1\\\kappa_\varphi^j > \eta}}^M \mu_Z\left[V^j \oplus (-Tx')\right]\right|\\
    & = \left|\sum\nolimits_{\substack{k = 1\\\kappa_\varphi^k > \eta}}^M \sum\nolimits_{\substack{j = 1\\\kappa_\varphi^j \leq \eta}}^M \mu_Z\left[\big(V^k \oplus (- Tx)\big) \cap \big(V^j \oplus (- Tx')\big)\right]\right.\\
    &\quad \left.- \sum\nolimits_{\substack{k = 1\\\kappa_\varphi^k \leq \eta}}^M \sum\nolimits_{\substack{j = 1\\\kappa_\varphi^j > \eta}}^M \mu_Z\left[\big(V^k \oplus (- Tx)\big) \cap \big(V^j \oplus (- Tx')\big)\right]\right|
\end{align*}
The additional restrictions below the sums exclude the case $k \neq j$, so that we can use the following considerations based on the notation introduced in the proof of Theorem~\ref{ThFLipschitz}:
\[
    \lim_{x' \to x} M^{k, j}(x, x')
    = M^{k, j}(x, x)
    = \begin{cases}
        0 & \text{if}\ k \neq j,\\
        \mu_Z\left[\bar{V}^k \oplus (- Tx)\right] & \text{if}\ k = j.
    \end{cases}
\]
We receive
\allowdisplaybreaks
\begin{align*}
    & \lim_{x' \to x} \left|\mathcal{Q}_{\mathrm{EP}_\eta}(x) - \mathcal{Q}_{\mathrm{EP}_\eta}(x')\right|\\
    & = \lim_{x' \to x} \left|\sum\nolimits_{\substack{k = 1\\\kappa_\varphi^k > \eta}}^M \sum\nolimits_{\substack{j = 1\\\kappa_\varphi^j \leq \eta}}^M \mu_Z\left[\big(V^k \oplus (- Tx)\big) \cap \big(V^j \oplus (- Tx')\big)\right]\right.\\
    &\quad \left.- \sum\nolimits_{\substack{k = 1\\\kappa_\varphi^k \leq \eta}}^M \sum\nolimits_{\substack{j = 1\\\kappa_\varphi^j > \eta}}^M \mu_Z\left[\big(V^k \oplus (- Tx)\big) \cap \big(V^j \oplus (- Tx')\big)\right]\right|\\
    & = \lim_{x' \to x} \left|\sum\nolimits_{\substack{k = 1\\\kappa_\varphi^k > \eta}}^M \sum\nolimits_{\substack{j = 1\\\kappa_\varphi^j \leq \eta}}^M M^{k, j}(x, x') - \sum\nolimits_{\substack{k = 1\\\kappa_\varphi^k \leq \eta}}^M \sum\nolimits_{\substack{j = 1\\\kappa_\varphi^j > \eta}}^M M^{k, j}(x, x')\right|
    = 0\tag*{$\qed$}
\end{align*}

\begin{proposition}
Assume \ref{A1}, \ref{A2} and let the support of $\mu_Z$ be compact. Fix $\alpha \in (0, 1)$, then the function 

$\mathcal{Q}_{\mathrm{VaR}_\alpha}(x) := \mathrm{VaR}_\alpha\left[\mathbb{F}(x)\right] = \inf\left\{\eta \in \mathbb{R} \; | \; \mu_Z\left[\left\{\omega \in \Omega \; | \; \mathbb{F}(x)(\omega) \leq \eta\right\}\right] \geq \alpha\right\}$

\noindent is continuous. Furthermore, let the induced feasible set $F_Z$ be nonempty and compact, then $\min_{x \in F_Z} \left\{\mathcal{Q}_{\mathrm{VaR}_\alpha}(x)\right\}$ is solvable.
\end{proposition}

{\it Proof} See the proof of \cite[Th. 2]{Iv14} in combination with Theorem~\ref{ThFLipschitz}. \qed

\section{Implications for qualitative stability} \label{Sec_stability}
We shall now examine the behaviour of local optimal values or optimal solutions sets of \eqref{StochasticBi-levelProgram} under perturbations of the underlying probability measure. As the measure's support may vary as well, stability analysis usually relies on the assumption of complete recourse, i.e. $\mathrm{dom}\; f = \mathbb{R}^n \times \mathbb{R}^s$ (cf. \cite{BuClDe20} for the case without integrality constraints). However, in the present setting, Lemma~\ref{LemmaDomPhi} and Lemma~\ref{LemmaA1DomainsPolyhedral} yield the existence of a polyhedral set $D \subsetneq \mathbb{R}^s$ such that $(x,z) \in \mathrm{dom} \, f$ holds if and only if $Tx+z \in D$. Thus, the domain of $f$ cannot encompass the whole space. Therefore, we shall only consider measures whose support is contained in the domain of function $f(x,\cdot)$ regardless of the leader's decision $x \in X$, i.e. confine the analysis to the set $\mathcal{M}(X) 
    := \left\{\mu \in \mathcal{P}(\mathbb{R}^s) \; | \; X \times \mathrm{supp}\; \mu \in \mathrm{dom}\; f\right\}$.
Since $f$ is measurable by Lemma~\ref{Lemma_phimeasurable2}, it follows $\left(\delta_x \otimes \mu\right) \circ f^{-1} \in \mathcal{P}(\mathbb{R})$, where $\delta_x \in \mathcal{P}(\mathbb{R}^n)$ denotes the Dirac measure at $x \in \mathbb{R}^n$. For fixed $x$, the mapping $f(x, \cdot)$ is bounded by Corollary~\ref{Cor_bounded}, so $\left(\delta_x \otimes \mu\right) \circ f^{-1}$ has a bounded support. Let the probability space $(\Omega, \mathcal{F}, \mathbb{P})$ be atomless (cf. \cite[Rem. 5.2]{BuClDe20}), then there exists some random variable $Z_{(x, \mu)} \in \mathcal{L}^p(\Omega, \mathcal{F}, \mathbb{P})$ such that $\left(\delta_x \otimes \mu\right) \circ f^{-1}
    = \mathbb{P} \circ Z_{(x, \mu)}^{-1}$.

Let $\rho: \mathcal{L}^p(\Omega, \mathcal{F}, \mathbb{P}) \to \mathbb{R}$ with $p \in [1, \infty)$ be a convex function which is nondecreasing and law-invariant. Thus, we may consider the mapping 

\mbox{$\mathcal{Q}_\rho: X \times \mathcal{M}(X) \to \mathbb{R}$} with $\mathcal{Q}_\rho(x, \mu)
    := \rho[Z_{(x, \mu)}]$,
    
\noindent which is well defined due to the law invariance of $\rho$.

In what follows, we will endow $\mathcal{P}(\mathbb{R}^s)$ with the \textit{topology of weak convergence}, i.e., the topology where a sequence $\{\mu_n\}_{n \in \mathbb{N}} \subset \mathcal{P}(\mathbb{R}^s)$ converges weakly to $\mu \in \mathcal{P}(\mathbb{R}^s)$ if and only if $    \lim_{n \to \infty}\int_{\mathbb{R}^s} h(t)~\mu_n(\mathrm{d}t)
    = \int_{\mathbb{R}^s} h(t)~\mu(\mathrm{d}t)$ is valid for any bounded continuous function $h: \mathbb{R}^s \to \mathbb{R}$.
    
A subset $\mathcal{M}
    \subseteq \mathcal{M}^p_s
    := \left\{\nu \in \mathcal{P}(\mathbb{R}^s) \; | \; \int_{\mathbb{R}^s} \|t\|^p~\mu(\mathrm{d}t) < \infty\right\}$, which denotes the set of Borel probability measures on $\mathbb{R}^s$ with finite moments of order $p \in [1, \infty)$, is called \textit{locally uniformly $\|\cdot\|^p$-integrating} if $\mu_n \to \mu \in \mathcal{M}$ implies $\lim_{n \to \infty}\int_{\mathbb{R}^s} \|t\|^p~\mu_n(\mathrm{d}t)
    = \int_{\mathbb{R}^s} \|t\|^p~\mu(\mathrm{d}t)$ for every sequence $\{\mu_n\}_{n \in \mathbb{N}} \subseteq \mathcal{M}$. Details and examples can be found in \cite[\S 5]{BuClDe20} and in \cite[\S 2]{ClKrSc17}.

The next two sets are adapted from the corresponding results in \cite{BuClDe20} and \cite{ClKrSc17}, respectively.
\begin{theorem}\label{Th_stability}
All assumptions as described above. Let $\mathcal{M} \subseteq \mathcal{M}(X)$ be locally uniformly $\|\cdot\|^p$-integrating and let $(x, \mu) \in X \times \mathcal{M}$ be such that $D_\varphi(x) - \{Tx\}$ is a $\mu-$null set. Then $\mathcal{Q}_\rho|_{X \times \mathcal{M}}$ is continuous at $(x, \mu)$ with respect to the product topology of the standard topology on $\mathbb{R}^n$ and the topology of weak convergence on $\mathcal{M}$.
\end{theorem}

{\it Proof} Based on Lemma~\ref{Lemma_phimeasurable2} and \eqref{nullset}, \cite[Cor. 2.3]{ClKrSc17} is applicable, where the required growth condition follows directly from Corollary~\ref{Cor_bounded} with $\gamma = 1$. \qed

\begin{remark}
If we assume alternatively all assumptions as described above, together with $\mathcal{M} \subseteq \mathcal{M}(X)$ is locally uniformly $\|\cdot\|^p$-integrating and $\mu \in \mathcal{M}$ is absolutely continuous with respect to the Lebesgue-Borel measure on $\mathbb{R}^s$ in the above theorem,  then $\mathcal{Q}_\rho|_{X \times \mathcal{M}}$ is continuous at $(x, \mu)$ for all $x \in X$ with respect to the product topology of the standard topology on $\mathbb{R}^n$ and the topology of weak convergence on $\mathcal{M}$.
\end{remark}

Next, we consider the parametric optimization problem
\begin{align}\label{Pmu}
    \min_x \left\{\mathcal{Q}_\rho(x, \mu) \; | \; x \in X\right\}.\tag{$P_\mu$}
\end{align}
As \eqref{Pmu} is non-convex, we are interested in the sets of local optimal solutions: For any open set $V \subseteq \mathbb{R}^n$, we introduce the extended real-valued \textit{localized optimal value function} $\xi_V: \mathcal{M}(X) \to \overline{\mathbb{R}}$ with 

$\xi_V(\mu) 
    := \inf_x \left\{\mathcal{Q}_\rho(x, \mu) \; | \; x \in X \cap \mathrm{cl}\; V\right\}$, 
    
\noindent and the \textit{localized optimal solution set mapping} $\Xi_V: \mathcal{M}(X) \rightrightarrows \mathbb{R}^n$ with

$\Xi_V(\mu)
    := \mathrm{Argmin}_x\left\{\mathcal{Q}_\rho(x, \mu) \; | \; x \in X \cap \mathrm{cl}\; V\right\}$.

\noindent Additional assumptions are required to study the stability of local optimal solution sets: The set $\Xi_V(\mu)$ is a so-called \textit{complete local minimizing (CLM) set} of \eqref{Pmu} with respect to $V$ if $\emptyset \neq \Xi_V(\mu) \subseteq V$ for a given $\mu \in \mathcal{M}(X)$ and an open set $V \subseteq \mathcal{R}^n$. 
\begin{theorem}\label{Th_stability2}
All assumptions as described above. Let $\mathcal{M} \subseteq \mathcal{M}(X)$ be locally uniformly $\|\cdot\|^p$-integrating and let $\mu$ has a uniformly bounded Lebesgue density. Then the following statements hold true:
\begin{enumerate}
    \item The restriction $\xi_{\mathbb{R}^n}|_\mathcal{M}$ is upper semicontinuous at $\mu \in \mathcal{M}$ with respect to the topology of weak convergence on $\mathcal{M}$.
\end{enumerate}
In addition, assume that $\mu \in \mathcal{M}$ is such that $\Xi_V(\mu)$ is a CLM set of \eqref{Pmu} with respect to some open bounded set $V \subset \mathbb{R}^n$. Then the following statements hold true:
\begin{enumerate}[start=2]
    \item The restriction $\xi_V|_\mathcal{M}$ is continuous at $\mu \in \mathcal{M}$ with respect to the topology of weak convergence on $\mathcal{M}$.
    \item The restriction $\Xi_V|_\mathcal{M}$ is upper semicontinuous at $\mu$ with respect to the topology of weak convergence on $\mathcal{M}$ in the sense of Berge (cf. \cite{Be59}), i.e. for any open set $O \subseteq \mathbb{R}^n$ with $\Xi_V(\mu) \subset O$ there exists a weakly open neighborhood $N$ of $\mu$ such that $\Xi_V(\mu_Z) \subseteq O$ holds for all $\mu_Z \in N \cap \mathcal{M}$.
\end{enumerate}
\end{theorem}

{\it Proof} Since the measure $\mu$ is endowed with a density, the assertions follow together with \cite[Cor. 2.4]{ClKrSc17}. \qed

\section{Conclusions}
We studied the function obtained by optimizing a linear functional over the set of minimizers of an integer linear problem and characterized its set of discontinuity points. This allowed us to formulate continuity results for the risk-averse bi-level stochastic model when the underlying Borel measure is absolutely continuous with respect to the Lebesgue measure. Moreover, we quantified the continuity by providing sufficient conditions for the Hölder continuity for a comprehensive class of risk measures. Qualitative stability results with respect to perturbations of the underlying probability measure conclude the paper.

Due to the boundedness of the function $\varphi$, it is easy to verify that all of our results for the optimistic model carry over to the pessimistic setting.

For risk measures defined on $L^1$, we have obtained an optimization problem with Lipschitz continuous goal function. Thus, Clarke's generalized derivatives can be used to formulate optimality conditions. However, the characterization of such derivatives was beyond the scope of the present paper and will be investigated in future work. Moreover, we will investigate algorithmic approaches based on these derivatives.

While the assumption that the feasible set of the lower level problem is finite arises natuarilly in a combinatorial or binary optimization setting, we plan to extend the analysis to unbounded models in future work.

\begin{acknowledgements}
The second author thanks the Deutsche Forschungsgemeinschaft for its support via the Collaborative Research Center TRR 154.
\end{acknowledgements}

\end{document}